\documentclass[preprint,amsmath,amssymb,aps,pre,nobibnotes,longbibliography]{revtex4-1}  

\usepackage{inputenc}
\usepackage{amsmath}
\usepackage{amssymb, amsfonts}
\usepackage{mathtools}
\usepackage{graphicx,epstopdf,bm}
\usepackage{tabularx}
\usepackage[colorlinks,bookmarks,urlcolor=blue,citecolor=blue,linkcolor=blue]{hyperref}
\usepackage[caption=false]{subfig}
\usepackage{algorithm}
\usepackage{algpseudocode}

\usepackage[T1]{fontenc}
\usepackage{bbm}
\usepackage{mathtools}
\usepackage{cases}
\usepackage{multirow}

\usepackage{bbm}

\newcommand{\diffophA}{\diffoph^\textrm{A}}
\newcommand{\diffophB}{\diffoph^\textrm{B}}
\newcommand{\diffophC}{\diffoph^\textrm{C}}

\newcommand{\phiA}{\phi^\textrm{A}}
\newcommand{\phiB}{\phi^\textrm{B}}
\newcommand{\phiC}{\phi^\textrm{C}}

\newcommand{\ZA}{Z_\textrm{A}}
\newcommand{\ZB}{Z_\textrm{B}}
\newcommand{\ZC}{Z_\textrm{C}}

\newcommand{\Lp}{\mathcal{L}}
\newcommand{\Lpb}{\mathcal{L}_{\textrm{b}}}
\newcommand{\kB}{k_\textrm{B}}

\newcommand{\kon}{k_\textrm{on}}
\newcommand{\koff}{k_\textrm{off}}
\newcommand{\RP}{R_\textrm{P}}
\newcommand{\RC}{R_\textrm{C}}
\newcommand{\piAB}{\pi_{\textrm{AB}}}
\newcommand{\piC}{\pi_{\textrm{C}}}
\newcommand{\pihatAB}{\hat{\pi}_{\textrm{AB}}}
\newcommand{\pihatC}{\hat{\pi}_{\textrm{C}}}

\newcommand{\ZhatA}{\hat{Z}_\textrm{A}}
\newcommand{\ZhatB}{\hat{Z}_\textrm{B}}
\newcommand{\ZhatC}{\hat{Z}_\textrm{C}}

\renewcommand{\vec}[1]{\boldsymbol{#1}}
\newcommand{\vecs}[1]{\boldsymbol{#1}}
\newcommand{\paren}[1]{\left(#1\right)}
\newcommand{\brac}[1]{\left[#1\right]}
\newcommand{\E}{\mathbb{E}}
\newcommand{\avg}[1]{\E[#1]}

\newcommand{\D}[2]{\frac{d#1}{d#2}}
\newcommand{\PD}[2]{\frac{\partial#1}{\partial#2}}

\newcommand{\abs}[1]{\left|#1\right|}
\newcommand{\norm}[1]{\Vert#1\Vert}

\DeclareMathOperator{\prob}{Pr}

\newcommand{\DA}{D^{\textrm{A}}}
\newcommand{\DB}{D^{\textrm{B}}}
\newcommand{\DC}{D^{\textrm{C}}}

\newcommand{\bi}{\vec{i}}
\newcommand{\bj}{\vec{j}}
\newcommand{\bk}{\vec{k}}

\newcommand{\ve}{\vec{e}}

\newcommand{\veta}{\vecs{\eta}}

\newcommand{\vp}{\vec{p}}
\newcommand{\vP}{\vec{P}}
\newcommand{\vC}{\vec{C}}

\newcommand{\va}{\vec{a}}
\newcommand{\vb}{\vec{b}}
\newcommand{\vc}{\vec{c}}
\newcommand{\vx}{\vec{x}}
\newcommand{\vy}{\vec{y}}
\newcommand{\vz}{\vec{z}}

\newcommand{\vf}{\vec{f}}
\newcommand{\vF}{\vec{F}}

\newcommand{\vv}{\vec{v}}
\newcommand{\vV}{\mathcal{V}}

\newcommand{\vO}{\vec{0}}

\newcommand{\vA}{\vec{A}}
\newcommand{\vB}{\vec{B}}
\newcommand{\rb}{r_{\textrm{b}}}

\newcommand{\vqa}{\vec{q}^{a}}
\newcommand{\vqb}{\vec{q}^{b}}
\newcommand{\vqc}{\vec{q}^{c}}

\newcommand{\vQa}{\vec{Q}^{a}}
\newcommand{\vQb}{\vec{Q}^{b}}
\newcommand{\vQc}{\vec{Q}^{c}}

\newcommand{\ind}{\mathbbm{1}}

\def\R{\mathbb{R}}

\newcommand{\pb}{p_{\textrm{b}}}
\newcommand{\pbO}{p_{\textrm{b},0}}

\newcommand{\Tb}{T_{\textrm{bind}}}

\newcommand{\equil}[1]{\bar{#1}}
\newcommand{\peq}{\equil{p}}
\newcommand{\pbeq}{\equil{p}_{\textrm{b}}}

\newcommand{\Dmat}{\mathcal{D}}
\newcommand{\Dbmat}{\mathcal{D}_{\textrm{b}}}
\newcommand{\Kd}{K_{\textrm{d}}}
\newcommand{\Id}{\mathcal{I}_d}

\renewcommand{\epsilon}{\varepsilon}
\renewcommand{\rb}{\epsilon}

\newcommand{\Pbound}{P_{\textrm{bound}}}

\newcommand{\kp}{\kappa^{+}}
\newcommand{\km}{\kappa^{-}}
\newcommand{\Rp}{\mathcal{R}^{+}}
\newcommand{\Rm}{\mathcal{R}^{-}}
\newcommand{\diffop}{\mathcal{L}}
\newcommand{\Rph}{\mathcal{R}^{+}_h}
\newcommand{\Rmh}{\mathcal{R}^{-}_h}
\newcommand{\diffoph}{\mathcal{L}_h}
\newcommand{\fabc}{f^{(a,b,c)}}
\newcommand{\Fabcijk}{F_{\bi^a \bj^b \bk^c}}
\newcommand{\vVabc}{\vV_{\bi^a \bj^b \bk^c}}

\newcommand{\kpijk}{\kappa_{i j k}^{+}}
\newcommand{\kpkij}{\kappa_{k\vert i j}^{+}}
\newcommand{\kpij}{\kappa_{i j}^{+}}
\newcommand{\kmk}{\kappa_{k}^{-}}
\newcommand{\kmijk}{\kappa_{i j k}^{-}}
\newcommand{\kmkij}{\kappa_{i j \vert k}^{-}}

\usepackage{amsthm}
\theoremstyle{plain}

\begin{document}

\title{An Unstructured Mesh Reaction-Drift-Diffusion Master Equation with Reversible Reactions}

\author{Samuel A. Isaacson}
\email{isaacson@math.bu.edu}
\affiliation{Department of Mathematics and Statistics, Boston University}

\author{Ying Zhang}
\email{ying1.zhang@northeastern.edu}
\affiliation{Department of Mathematics and Department of Biology, Northeastern University}

\numberwithin{equation}{section}


\begin{abstract}
	We develop a convergent reaction-drift-diffusion master equation (CRDDME) to facilitate the study of reaction processes in which spatial transport is influenced by drift due to one-body potential fields within general domain geometries. The generalized CRDDME is obtained through two steps. We first derive an unstructured grid jump process approximation for reversible diffusions, enabling the simulation of drift-diffusion processes where the drift arises due to a conservative field that biases particle motion. Leveraging the Edge-Averaged Finite Element method, our approach preserves detailed balance of drift-diffusion fluxes at equilibrium, and preserves an equilibrium Gibbs-Boltzmann distribution for particles undergoing drift-diffusion on the unstructured mesh. We next formulate a spatially-continuous volume reactivity particle-based reaction-drift-diffusion model for reversible reactions of the form $\textrm{A} + \textrm{B} \leftrightarrow \textrm{C}$. A finite volume discretization is used to generate jump process approximations to reaction terms in this model. The discretization is developed to ensure the combined reaction-drift-diffusion jump process approximation is consistent with detailed balance of reaction fluxes holding at equilibrium, along with supporting a discrete version of the continuous equilibrium state. The new CRDDME model represents a continuous-time discrete-space jump process approximation to the underlying volume reactivity model. We demonstrate the convergence and accuracy of the new CRDDME through a number of numerical examples, and illustrate its use on an idealized model for membrane protein receptor dynamics in T cell signaling.
\end{abstract}

\maketitle

\section{Introduction}
\label{sec:intro}
Mathematical and computational models have become an invaluable tool in understanding the dynamics of many stochastic cellular processes, governed by a co-action between spatial transport and chemical reactions. Examples of such processes include how ligand binding can affect the functioning of membrane-bound proteins \cite{SuEtAl2022,AlbertsMOLECELLBIO} and how tethered catalytic reactions control the initiation and integration of T cell activation \cite{Zhang2019}. One challenge in developing mathematical models of such processes is that spatial transport of molecules occurs within heterogeneous environments, such as cell membranes, which can significantly impact the dynamics of chemical reaction processes \cite{LarsenEtAl2020,RojasEtAl2021}. For example, the geometry of cellular reaction domains can dramatically alter the initial phase of cell polarization \cite{GieseEtAl2015}, can influence gradient sensing of migratory cells \cite{ChengEtAl2020}, and can modulate the propagation of signals from the cell membrane to nucleus \cite{IsaacsonPLOS20}.

An additional modeling challenge arises from experimental observations that in many biochemical systems the spatial dynamics of diffusing particles also involves drift. At the single cell scale, sources of drift include heterogeneity in the cytoplasm, interactions with cellular structures, variations in chemical potential, and active transport. For instance, during T cell synapse formation protein diffusion is influenced by both membrane geometry and the concentration of F-actin fibers near the cell membrane \cite{DemondEtAl2008,BasuHuse2017,Siokisetal2018}. The latter can be modeled as imparting drift due to a conservative force, i.e. a background/one-body potential field, to proteins within the cell membrane. One body background potentials have also played critical roles in developing models for molecular-motor based active transport \cite{PeskinOsterMotor1995,AtzbergerPeskin2006}, in modeling protein movement subject to the effect of volume exclusion by chromatin \cite{IsaacsonPNAS2011}, and in approximating reflecting boundary conditions within cells to prevent proteins from crossing domain boundaries \cite{AtzbergerIsaacson2009}. To facilitate the study of reaction-drift-diffusion processes with drift due to background potentials, it is necessary to develop microscopic models and numerical simulation methods that can accurately and explicitly resolve different physical mechanisms modulating the movement of molecules within realistic cellular domains.

In this work, we continue with the setup introduced in \cite{ZhangPhDthesis,IsaacsonCRDME2013,IsaacsonZhang17} and focus on the spatially continuous volume reactivity (VR) model \cite{Flegg:2013bc,Schoneberg:2014jr,PrustelMS2014,ErbanChapman,Donevetal2018,Zhang2019,StrobergSchnell2021,HuhnEtAl2024}, modeling particle motion as a drift-diffusion process where the drift is due to a one-body background potential field (i.e. in the absence of reaction, particles move independently with statistics governed by the Fokker-Planck equation). A bimolecular association reaction of the form $\textrm{A} + \textrm{B} \to \textrm{C}$ is modeled through a binding kernel that determines the probability density per time an \textrm{A} molecule at $\vx$ and a \textrm{B} molecule at $\vy$ react and produce a \textrm{C} molecule at $\vz$ inside the domain. 
First order reactions of the form $\textrm{A} \to \textrm{B}$ are treated as an internal process, with individual particles switching type based on exponential clocks. A special case of the general VR model, for a specific choice of the bimolecular association kernel, is the Doi model \cite{TeramotoDoiModel1967,DoiSecondQuantA,DoiSecondQuantB}.

In the VR model, a coupled system of stochastic processes represents the time-evolution of the state of the chemical system, i.e. of the number of particles of each chemical species and their respective positions. These processes can be described via partial integral differential equations (PIDEs) for the probability density they have a given value at a given time, i.e. forward Kolmogorov equations associated with the processes which Doi presented via his Fock space representation~\cite{DoiSecondQuantA,DoiSecondQuantB}. Alternatively, weak representations for the dynamics of the processes can be formulated as measure-valued processes~\cite{IsaacsonMaMeanfield2022,HeldmanMeanFieldDrift2024}.

When considering such models for systems that are closed and only involve reversible reactions, such as $\textrm{A} + \textrm{B} \leftrightarrows \textrm{C}$, from statistical mechanics we expect detailed balance of reaction fluxes to hold at equilibrium~\cite{Donevetal2018,FrohnerNoe2018,ZhangIsaacson2022}. This in turn imposes a functional relationship between the association and dissociation reaction kernels~\cite{ZhangIsaacson2022,FrohnerNoe2018}, enabling one to determine dissociation reaction kernels in terms of assocation kernels (or vice-versa). If this relationship is not preserved in numerical approximations to the VR model, one can obtain inaccurate approximations to the underlying dynamics and statistics~\cite{ZhangIsaacson2022,FrohnerNoe2018}.

One common approach for generating approximate realizations of the stochastic process of diffusing and reacting particles in the VR model is through Brownian Dynamics methods that discretize these processes in time \cite{AndrewsBrayPhysBio2004,ReaddyPLOS2013,ReaddyPLOS2019}. To account for interacting particle dynamics, a Monte-Carlo based algorithm was proposed in \cite{FrohnerNoe2018}. This scheme includes a Metropolis-Hastings sampling step that ensures the resulting dynamics satisfy detailed balance of reaction fluxes for interacting particles. While this ensures detailed balance of reaction fluxes, the stochastic differential equation solvers that propagate particles in space would not be expected to preserve detailed balance of drift-diffusion fluxes or the Gibbs-Boltzmann distribution at equilibrium. More recently, for surface reaction-drift-diffusion dynamics a novel method that discretizes in space and time was proposed in~\cite{AtzbergerDriftDiffusion2022}. This latter method is based on Smoluchowski reaction dynamics, and was constructed to preserve detailed balance of drift-diffusion fluxes within the discretized model. 

\begin{figure}[!tbp]
	\centering
	\includegraphics[width = 0.99\textwidth]{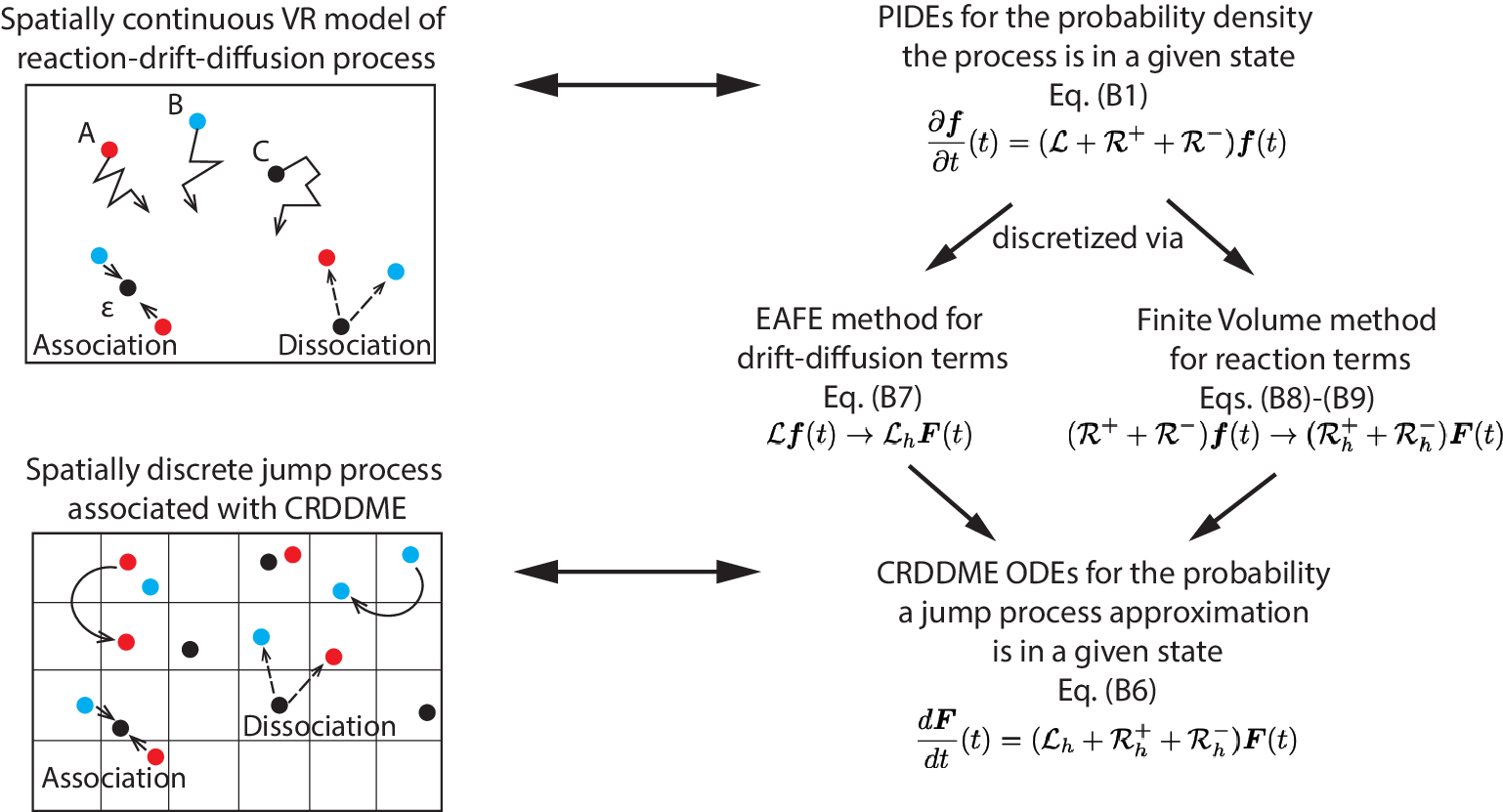}
	\caption{Our approach to approximating the stochastic reaction-drift-diffusion process of the Volume Reactivity (VR) model (upper left) by a spatially-discrete jump process (lower left). Derivation proceeds clockwise from the former to the latter via the CRDDME we develop in this work.}
	\label{fig:CRDDME_flowchart}
  \end{figure}  
In this work we take a different approach to approximate the stochastic processes associated with the VR model via developing a convergent spatially-discrete jump process approximation. We begin by developing an unstructured grid spatial discretization of the system of PIDEs that describe the forward Kolmogorov equation for the VR model. Our discretization is formulated so that the resulting set of ODEs we obtain corresponds to a master equation, i.e. the forward Kolmogorov equation for a jump process. We call this master equation the Convergent Reaction-Drift-Diffusion Master Equation (CRDDME). The CRDDME then represents a new, spatially discrete, stochastic reaction-drift-diffusion model that approximates the underlying spatially continuous VR model. In the jump process associated with the CRDDME, particles  move by hopping between mesh sites with fixed probabilities per time, and can react with probabilities per time based on their positions and separation within the discrete mesh (for bimolecular reactions). Due to its high dimensionality, for all but the smallest systems it is not possible to directly solve the ODEs of the CRDDME numerically. Instead, our final numerical methods use Stochastic Simulation Algorithms (SSAs), i.e. Gillespie Methods, to sample exact realizations of the associated jump processes. This overall approach to derive an approximation to the spatially-continuous VR model process is illustrated in Fig.~\ref{fig:CRDDME_flowchart}.

In this way we obtain a convergent method that can be used to simulate stochastic reaction-drift-diffusion processes in complex 2D or 3D cellular geometries defined by unstructured meshes, as often reconstructed from experimental imaging assays. The resulting CRDDME model supports variable diffusivity and drift due to potential fields, converges to the PIDEs of the corresponding VR model as the mesh spacing is taken to zero, and is consistent with the popular reaction-diffusion master equation (RDME) in its handling of linear reactions and spatial transport. It represents a generalization of the convergent RDME (CRDME) we previously developed in \cite{IsaacsonCRDME2013,IsaacsonZhang17} to include drift due to one-body potential fields and spatially varying diffusivities. The CRDDME also avoids an implicit challenge in using RDME models in two or more dimensions, where bimolecular reactions are lost as the mesh spacing is taken to zero~\cite{IsaacsonRDMELims,Hellander:2012jk,IsaacsonCRDME2013}. In using RDME models, one must then search, and hope, for the existence of mesh scales that are not so small that bimolecular reactions are inaccurately resolved (i.e. lost), but not so large that spatial transport is inaccurately approximated. In contrast, when using the CRDDME one can simply choose the mesh spacing to achieve a desired level of accuracy in approximating the underlying VR model. The CRDDME is also constructed so that in closed systems, it is consistent with equilibrium properties of the underlying particle model, including detailed balance of drift-diffusion fluxes, an equilibrium Gibbs-Boltzmann distribution in the absence of reactions, and detailed balance of reversible reaction fluxes. 

To construct a CRDDME that incorporates reversible reactions on an unstructured mesh, we utilize a hybrid finite element / finite volume discretization approach as illustrated in Fig.~\ref{fig:CRDDME_flowchart}. We begin in the next section by approximating the Fokker-Planck equation for a single particle within a bounded domain by a spatial jump process, adapting the Edge-Average Finite Element (EAFE) Method of \cite{XuZikatanov1999} to derive spatial transition rates for one particle to hop between neighboring mesh voxels. Here the particle hops between voxels of the polygonal (barycentric) dual mesh to a triangulated/tetrahedral mesh of the domain. The resulting discretized drift-diffusion operator defines a proper transition rate matrix under the same Delaunay-type mesh conditions as needed in the pure diffusion case \cite{LotstedtFERDME2009}, while also reducing to the transition rate matrix derived for the pure diffusion case in \cite{LotstedtFERDME2009} when there is no potential field. Our final discretization can be interpreted as an unstructured mesh analog to the Cartesian grid finite difference approximations for drift-diffusion processes arising from potential fields developed in \cite{ElstonPeskinJTB2003,IsaacsonPNAS2011}, and Cartesian-grid finite volume approximation developed in \cite{Schutteetal2011}. We demonstrate in Appendix~\ref{app:EAFEAccuracy} the accuracy of the EAFE discretization by considering several steady-state problems with different choices of potential functions. As a partial differential equation (PDE) discretization, the EAFE method is shown to be second order accurate.

In Section \ref{sec:RevRX} we then approximate the reversible $\textrm{A} + \textrm{B} \leftrightarrows \textrm{C}$ bimolecular reaction on the same underlying polygonal mesh. We begin by introducing the abstract continuous VR model for reversible reactions in Section \ref{sec:absContModel}. In Section \ref{sec:DBContModel} we show how association and dissociation kernels should be related to be consistent with detailed balance, using this relationship to define dissociation kernels in terms of association kernels. We apply a finite volume method to derive a jump process approximation to the association and dissociation reaction terms in the spatially-continuous VR model in Section \ref{sec:FVM_RXs}. Combining the finite element discretization of the drift-diffusion terms from Section \ref{sec:DriftDiffApprox} and the finite volume approximation of reaction terms from Section \ref{sec:RevRX}, we summarize in Section \ref{sec:discreteModel} a complete spatially-discrete CRDDME model for a pair of particles that undergoes drift-diffusion with the reversible reaction. We show in Appendix \ref{app:multipartModel} how this can be generalized to a multi-particle system through formulating the general multiparticle VR model with drift, and discretizing this model to obtain a multiparticle CRDDME. In particular, the resulting multiparticle transition rates arise from the two-particle rates in a simple manner. In Appendix \ref{app:multipartDB} we show these multiparticle transition rates are consistent with detailed balance of reaction fluxes within general multiparticle models (when the two-particle transition rates are also chosen to satisfy detailed balance). Finally, in Section \ref{sec:numerical_examples} we consider a number of numerical examples to demonstrate the convergence and accuracy of the CRDDME in approximating the continuous VR model, and to illustrate how the CRDDME can be applied to biological systems, by considering a model for TCR-pMHC dynamics in early stage immune synapse formation that incorporates key components from several recent studies~\cite{DharanFarago2017,Siokis2017,Siokisetal2018}. 

We note that in the remainder, convergence refers to convergence of the solution to the CRDDME for the probability the jump process approximation is in a given state to the solution of the PIDEs for the probability density the VR model stochastic process is in a given state, or to convergence of statistics of the CRDDME model to those of the VR model. This represents a type of weak convergence at the level of the jump process associated with the CRDME to the process associated with the VR model; we do not consider stronger types of convergence at the level of the associated processes in this work (such as convergence in probability).

Finally, note that while in this work we illustrate our approach for deriving a CRDDME on the reversible $\textrm{A} + \textrm{B} \leftrightarrows \textrm{C}$ reaction, our approach is applicable to general non-equilibrium systems of zero, first, and second order reactions. We choose to focus on the former case for simplicity of exposition, and to show the CRDDME can be formulated to respect the additional constraint of detailed balance that holds for such systems.

\section{Jump Process Approximation of the Fokker-Planck Equation with One-body Potential}
\label{sec:DriftDiffApprox}
We begin by deriving an unstructured mesh jump process approximation to the drift-diffusion of individual particles. Spatial transition rates between lattice sites are obtained by using an EAFE Method discretization of the Fokker-Planck operator on triangulated meshes \cite{XuZikatanov1999}. In this section, we summarize the method and derive the associated stiffness matrix that generates the transition rate matrix for the jump process. Readers interested in a more detailed discussion of the EAFE method should see \cite{XuZikatanov1999}.

In the absence of chemical reactions, the spatial movements of individual particles are independent processes that do not affect each other; and therefore, to obtain the spatial transition rates for multi-particle systems, it is sufficient to derive an approximation for a system with only one molecule that moves within a closed and bounded domain $\Omega \subset \R^d$, $d \leq 3$ \cite{IsaacsonSJSC2006}.
To formulate the underlying equations describing the state of the system, we denote by $p(\vx,t)$ the probability density the diffusing molecule is located at $\vx$ at time $t$ with initial condition $p(\vx,0) = p_0(\vx)$. The diffusivity of the moving molecule is denoted by $D(\vx)$, where
\begin{align*}
  \inf_{\vx \in \Omega} D(\vx) > C > 0,
\end{align*}
i.e. $D(\vx)$ is bounded from below by a positive constant, $C$, on $\Omega$. The background potential field that the molecule experiences while diffusing is given by $\phi(\vx)$. Here we assume that the potential field only varies in space (i.e. independent of time).
With these definitions, let $\diffop$ denote the drift-diffusion operator
\begin{equation}
  \diffop p(\vx,t) = \nabla\cdot\brac{D(\vx)\paren{\nabla p(\vx,t)+p(\vx, t)\nabla\phi(\vx)}}.
  \label{eq:FPoperator}
\end{equation}
Here and in the remainder we have assumed the potential is non-dimensional, absorbing a typical factor of $(\kB T)^{-1}$ into it, along with assuming the transport model is such that the diffusivity multiplies both the drift and diffusion terms. Such a formulation arises in a variety of physical applications, and can be derived by assuming a local Einstein relation holds at each point in space~\cite{GuerinDean2015}. It also reduces to the typical formulation that arises in the case of constant diffusivity~\cite{ElstonPeskinJTB2003}, which we will shortly restrict to when adding in reaction processes. With $\diffop$ defined, $p(\vx,t)$ then satisfies
\begin{equation}
	\begin{aligned}
  		\PD{p}{t}(\vx,t) &= \diffop p(\vx,t), &\vx \in \Omega \\
  		D(\vx)\paren{\nabla p(\vx, t)+p(\vx, t)\nabla\phi(\vx)}\cdot\veta(\vx) &= 0, &\vx \in \partial\Omega.
  	\end{aligned}
  	\label{eq:FPeqn}
\end{equation}
Here $\veta(\vx)$ denotes the outward normal to a surface at the point $\vx$.

We begin by showing how a standard piecewise linear finite element method (FEM) approximation to $\diffop$ in \eqref{eq:FPoperator} gives a discretization that can not be interpreted as corresponding to a transition rate matrix associated with a jump process. We denote by $\Omega_h$ a triangulated mesh approximation to $\Omega$, with $\mathcal{T}_h$ the collection of $K$ triangles/tetrahedra within $\Omega_h$, and let $\{\psi_i\}_{i=1}^N$ label the standard set of $N$ piecewise linear finite element method basis functions on $\Omega_h$. For each node/vertex of the mesh, labeled by $\{\vx_i\}_{i=1}^K$, we have that
\begin{equation*}
	\psi_j(\vx_i) = \delta_{i,j}.
\end{equation*}
We approximate $p(\vx,t)$ by
\begin{equation}
	p_h(\vx,t) = \sum_{i = 1}^Np_i(t)\psi_i(\vx),
	\label{eq:FEMapprox}
\end{equation}
and let $\vp_h(t) = \brac{p_h(\vx_1,t), \dots, p_h(\vx_N,t)}^T$ denote the vector of nodal values. A standard piecewise linear finite element discretization of \eqref{eq:FPeqn} then gives
\begin{equation}
	\D{\vp_h}{t} = M^{-1}S\vp_h,
	\label{eq:FEMrho}
\end{equation}
Here the mass matrix, $M$, is defined by
\begin{equation}
	M_{ij} = \int_{\Omega}\psi_i(\vx)\psi_j(\vx) \, d\vx,
	\label{eq:FPMass}
\end{equation}
and the stiffness matrix, $S$, is defined by
\begin{equation}
	S_{ij} = -\int_{\Omega}\brac{D(\vx)\paren{\nabla \psi_i(\vx)+\psi_i(\vx)\nabla\phi(\vx)}}\cdot \nabla \psi_j \, d\vx.
	\label{eq:FPStiffness}
\end{equation}
We note that the stiffness matrix defined above is not guaranteed to satisfy $S_{ij} \geq 0$ for $i \neq j$, even for a Delaunay triangulated mesh. As such, we can not generate a transition rate matrix from it to represent the rates at which the particle hops between mesh nodes (i.e. some transition rates could end up with negative values). In other words, we can not interpret the preceding semi-discretized model as corresponding to an underlying jump process for the particle hopping on the mesh.

To resolve this issue, we derive our stiffness matrix using the EAFE scheme developed in \cite{XuZikatanov1999}. Let $\omega_{ij}$ denote the entries of the standard piecewise linear finite element stiffness matrix for the Laplacian (i.e. without background potential field). Label by $E_{ij}$ the triangle edge connecting $\vx_j$ to $\vx_i$. The EAFE stiffness matrix is then given by
\begin{equation}
	S_{ij} = \begin{cases}
			\omega_{ij}\brac{\frac{1}{\abs{E_{ij}}}\int_{E_{ij}}\frac{e^{\phi(\vx(s))-\phi_j}}{D(\vx(s))} \, ds}^{-1}, &i \neq j \\
			-\sum_{k\neq j}S_{kj}, &i = j,
	\end{cases}
	\label{eq:EAFE_stiffMat}
\end{equation}
where $\abs{E_{ij}}$ denotes the length of the edge and $\phi_i = \phi(\vx_i)$. Note, in the case of a constant potential this reduces to the standard piecewise linear finite element stiffness matrix approximation to the Laplacian.

The system \eqref{eq:FEMrho} can be further simplified by using numerical quadrature to evaluate the elements of the mass matrix. One common choice, equivalent to mass lumping as previously used in the purely diffusive case \cite{LotstedtFERDME2009} is to use the trapezoidal rule, i.e.
\begin{equation}
	\int_{\Omega}\psi_i\psi_j \, d\vx 
		\approx \sum_{K\in\mathcal{T}_h}\frac{\abs{T}}{(d+1)}\sum_{m = 1}^{d+1}\delta_{\vx_i\vx_K^m}\delta_{\vx_j\vx_K^m},
\label{eq:lumpQuad}
\end{equation}
where $\vx_K^m$ denotes the $m$th node of the $K$th triangle. Using this approximation, we replace $M$ with a lumped mass matrix, $\Lambda$, with entries
\begin{equation}
	M \approx \Lambda_{ij} = \begin{cases}
		0, & i\neq j\\
		\abs{V_i}, &i = j,
	\end{cases}
	\label{eq:newQuadME}
\end{equation}
where in two-dimensions, $\abs{V_i}$ gives the area, and in three-dimensions the volume, of the $i$th (barycentric) dual mesh element. Inverting $\Lambda$ we obtain a spatially-discrete approximation to the Fokker-Planck equation
\begin{equation}
	\D{\vp_h}{t} = \Lambda^{-1}S \vp_h.
	\label{eq:semidiscreteFPE}
\end{equation}
We approximate by $\vP(t) = \brac{P_1(t), \dots, P_N(t)}^T = \Lambda\vp_h(t)$ the probabilities to be in each mesh element at time $t$. \eqref{eq:semidiscreteFPE} then simplifies to
\begin{equation}
	\D{\vP}{t} = \diffop_h\vP,
	\label{eq:semidisFEMProb}
\end{equation}
where we subsequently called $\diffop_h$ the discrete drift-diffusion operator given by
\begin{equation}
	\diffop_h := S\Lambda^{-1}.
	\label{eq:discreteDDOperator}
\end{equation}
Here the subscript $h$ indicates this is the spatially-discretized operator. Later in the section, $h$ will also represent the maximum diameter among all mesh voxels.

In the purely diffusive case, the operator $\diffop_h$ was shown to correspond to the transition rate matrix associated with a jump process whenever the triangulation is Delaunay in two dimensions, or under a condition where the dihedral angle at an edge between two facets of any tetrahedron in the mesh is non-obtuse in three dimensions \cite{LotstedtFERDME2009}. This is equivalent to the pure-diffusion stiffness matrix being an M-matrix \cite{XuZikatanov1999}. Our drift-diffusion $\diffop_h$ in~\eqref{eq:discreteDDOperator} is also guaranteed to represent a transition rate matrix when $S$ corresponds to an M-matrix. As such, our discretized Fokker-Planck operator, $\diffop_h$, will correspond to a transition rate matrix on any mesh for which the purely diffusive operator from \cite{LotstedtFERDME2009} represents a transition rate matrix.

When $\diffop_h$ corresponds to a transition rate matrix, for $i \neq j$ and $S_{i j}$ given by the EAFE discretization~\eqref{eq:EAFE_stiffMat}, we have that
\begin{equation}
	(\diffop_h)_{ij} = \frac{S_{ij}}{\abs{V_i}}
	\label{eq:transitionRate}
\end{equation}
can be interpreted as the probability per time (i.e. transition rate) for a molecule in the $j$th mesh element to hop to the $i$th mesh element, while
\begin{equation}
	-(\diffop_h)_{ii} = \frac{1}{\abs{V_i}}\sum_{j\neq i}S_{ji}
\end{equation}
gives the total probability per time for a molecule in $V_i$ to hop to any neighboring mesh element.

\subsection{Equilibrium of the Discretized Fokker-Planck Equation}
The equilibrium solution to the single-particle Fokker-Planck equation \eqref{eq:FPeqn} is given by the Gibbs-Boltzmann distribution
\begin{equation*}
	\bar{p}(\vx) = \frac{e^{-\phi(\vx)}}{\int_{\Omega}e^{-\phi(\vx)}\, d\vx}.
\end{equation*}
We now show that an analogous discrete probability distribution is an equilibrium solution for the EAFE-based jump process master equation \eqref{eq:semidisFEMProb}.

From statistical mechanics / thermodynamics, we expect that the equilbrium for our spatial discrete model~\eqref{eq:semidisFEMProb} should be given by a discrete Gibbs-Boltzmann probability distribution. Let $\bar{\vP}_i$ denote the probability the particle is in the $i$th voxel at equilibrium, with $\phi_i := \phi(x_i)$ the (discrete) potential in the voxel. The discrete Gibbs-Boltzmann distribution is then 
 \begin{equation}
	\bar{\vP}_i = \frac{e^{-\phi_i}\abs{V_i}}{Z},
\end{equation}
where $Z$ denotes the partition function 
\begin{equation*}
	Z = \sum_j e^{-\phi_j}\abs{V_j}.
\end{equation*}
Let $\bar{\vP} = \paren{\bar{\vP}_1, \bar{\vP}_2, ..., \bar{\vP}_N}$.

Detailed balance of the spatial transition rates holds for this equilibrium state, i.e.
\begin{equation*}
	\begin{aligned}
		(\diffop_h)_{ij}\bar{\vP}_j &= \frac{S_{ij}}{\abs{V_j}}\frac{e^{-\phi_j}\abs{V_j}}{Z} \\
																&= \frac{\omega_{ij}}{Z}\brac{\frac{1}{\abs{E_{ij}}}\int_{E_{ij}}\frac{e^{\phi(\vx(s))}}{D(\vx(s))} \, ds}^{-1} \\
																&= \frac{\omega_{ji}}{Z}\brac{\frac{1}{\abs{E_{ji}}}\int_{E_{ji}}\frac{e^{\phi(\vx(s))}}{D(\vx(s))} \, ds}^{-1} \\
																&= \frac{S_{ji}}{\abs{V_i}}\frac{e^{-\phi_i}\abs{V_i}}{Z} \\
																&= (\diffop_h)_{ji}\bar{\vP_i},
	\end{aligned}
\end{equation*}
where we have used that the pure-diffusion stiffness matrix is symmetric (i.e. $S_{ij} = S_{ji}$). Likewise, the Gibbs-Boltzmann distribution is a steady-state for the master equation model, i.e.
\begin{equation*}
	\begin{aligned}
		(\diffop_h\bar{\vP})_i &= \frac{1}{Z}\sum_j\frac{S_{ij}}{\abs{V_j}}\abs{V_j}e^{-\phi_j} \\
													 &= \frac{1}{Z}\sum_{j\neq i}\brac{\frac{\omega_{ij}}{\abs{E_{ij}}}\int_{E_{ij}}\frac{e^{\phi(\vx(s))}}{D(\vx(s))} \, ds}^{-1} - \frac{1}{Z}\sum_{j\neq i}\brac{\frac{\omega_{ji}}{\abs{E_{ji}}}\int_{E_{ji}}\frac{e^{\phi(\vx(s))}}{D(\vx(s))} \, ds}^{-1} \\
													 &= 0.
	\end{aligned}
\end{equation*}
We therefore observe that our semi-discrete jump process master equation model is consistent with the equilibrium solution to the continuous single-particle Fokker-Planck equation model, with detailed balance of the spatial fluxes holding in that equilibrium state. In Appendix~\ref{app:EAFEAccuracy} we demonstrate the accuracy of the EAFE method in the $L^2$ norm as a PDE discretization in the absence of sampling error by considering several steady state PDE problems on a square and a circle with different choices of potential functions. We observe that the rate of convergence of the EAFE method solution is second order as the mesh spacing is decreased.

\section{Reversible Reactions on Unstructured Meshes}
\label{sec:RevRX}
With the jump process approximation to the spatial movement of particles established, we now handle chemical reactions using the finite-volume-based jump process approximation approach for reversible bimolecular reactions we developed in \cite{IsaacsonZhang17}. For simplicity of notation, we derive a jump process approximation for the reversible bimolecular reaction $\textrm{A} + \textrm{B} \rightleftharpoons \textrm{C}$ in a system with only one $\textrm{A}$ and one $\textrm{B}$ molecule, or one $\textrm{C}$ complex, where each molecule experiences a background potential $\phiA(\vx)$, $\phiB(\vy)$, and $\phiC(\vz)$ respectively. In Appendix \ref{app:multipartModel} we show that the two-particle system can be generalized to a multi-particle system involving arbitrary numbers of $\textrm{A}$, $\textrm{B}$, and $\textrm{C}$ molecules, with the multi-particle reactive transition rates arising from the two-particle rates. This is also illustrated in Table \ref{Tab:CRDDMETransitionRates}, where it is apparent that the multiparticle reactive transition rates arise from the two-particle reactive transition rates multiplied by the number of distinct combinations of substrates. We begin in the next section by formulating the spatially-continuous two-particle reaction-drift-diffusion VR model. We then illustrate in Section \ref{sec:DBContModel} how reactive interaction functions for the reversible reaction should be related in the continuous model to be consistent with detailed balance of reaction fluxes holding at equilibrium. In Section \ref{sec:FVM_RXs} we use the finite volume method of \cite{IsaacsonZhang17} to obtain a master-equation-based discretization of the reaction terms in the continuous VR model, with the full hybrid EAFE-finite-volume based reaction-diffusion master equation approximation presented in Section \ref{sec:discreteModel}. Finally, we demonstrate in Section \ref{sec:DBCRDDME} that rates given by direct discretization do not satisfy detailed balance after discretization and instead, one should choose the reaction rates by enforcing detailed balance of reaction fluxes at the discrete level. By doing so, we obtain rates that are different than those by direct discretization, but can be interpreted as a quadrature approximation to them that preserves the rate of convergence of the full CRDDME approximation as the mesh width is taken to zero.

\subsection{Abstract Continuous Model}
\label{sec:absContModel}
We begin by considering the reversible bimolecular reaction $\textrm{A} + \textrm{B} \to \textrm{C}$ in a system whose state is either one $\textrm{A}$ and one $\textrm{B}$ molecule or one $\textrm{C}$ complex. As in Section \ref{sec:DriftDiffApprox}, we assume the reaction takes place within a $d$-dimensional bounded domain $\Omega \subset \R^d$, $d\leq 3$. Let $\vx\in\Omega$ denote the position of the $\textrm{A}$ molecule, $\vy\in\Omega$ the position of the $\textrm{B}$ molecule, and $\vz\in\Omega$ the position of the $\textrm{C}$ molecule. To formulate the underlying equations describing the state of the system, we denote by $p(\vx,\vy,t)$ the probability density the \textrm{A} and \textrm{B} molecules are unbound and located at $\vx$ and $\vy$ respectively at time $t$, and $\pb(\vz,t)$ the probability density the molecules are bound and the corresponding \textrm{C} molecule is located at $\vz$ at time $t$. For simplicity, the diffusivities of the molecules are assumed to be constant and are given by $\DA(\vx) = \DA$, $\DB(\vy) = \DB$, and $\DC(\vz) = \DC$. With $\Id$ the $d-$dimensional identity matrix, we define two constant diffusivity matrices, given by the block matrices
\begin{equation}
  \label{eq:diffusivityDefs}
  \begin{aligned}
    \Dmat &:= \begin{bmatrix}
      \DA  \Id & \vO \\
      \vO & \DB \Id
    \end{bmatrix} \\
    \Dbmat &:= \DC \Id.
  \end{aligned}
\end{equation}
Finally, with $\nabla_{\vx,\vy} = (\nabla_{\vx},\nabla_{\vy})$, $\nabla\Phi_{\vx,\vy} = (\nabla_{\vx}\phiA(\vx),\nabla_{\vy}\phiB(\vy))$, and $\nabla\Phi_{\vz} = \nabla_{\vz}\phiC(\vz)$, spatial operators in $(\vx,\vy)$ and $\vz$ are given by
\begin{equation}
	\begin{aligned}
		\Lp &= \nabla_{\vx,\vy}\cdot\Dmat\brac{\nabla_{\vx,\vy} + \nabla\Phi_{\vx,\vy}} \\
		\Lpb &= \nabla_{\vz}\cdot\Dbmat\brac{\nabla_{\vz}+\nabla\Phi_{\vz}}.
	\end{aligned}
	\label{eq:spatialOpsRevRx}
\end{equation}

The forward association $\textrm{A} + \textrm{B} \to \textrm{C}$ reaction process is defined by the probability density per unit time a reaction occurs creating a \textrm{C} molecule at $\vz$ given an \textrm{A} molecule at $\vx$ and a \textrm{B} molecule at $\vy$, denoted by $\kp(\vz\vert\vx,\vy)$. As such, the probability per time an $\textrm{A}$ molecule at $\vx$ and a $\textrm{B}$ molecule at $\vy$ successfully react to produce a $\textrm{C}$ molecule within the domain is
\begin{equation}
 \kp(\vx, \vy) := \int_{\Omega}\kp(\vz\vert\vx,\vy) \, d\vz.
\end{equation}

Similarly, we denote by $\km(\vx,\vy\vert\vz)$ the probability density per unit time a dissociation reaction occurs successfully and creates an $\textrm{A}$ molecule at $\vx$ and a $\textrm{B}$ molecule at $\vy$ given a $\textrm{C}$ molecule at $\vz$. With this definition, we let $\km(\vz)$ be the probability per time a \textrm{C} molecule at $\vz$ unbinds and produces \textrm{A} and \textrm{B} molecules within $\Omega$
\begin{equation}
	\km(\vz) := \int_{\Omega^2} \km(\vx,\vy\vert\vz) \, d\vx \, d\vy,
\end{equation}
where $\Omega^2 = \Omega \times \Omega \subset \R^{2d}$.

Given the preceding definitions, we can formulate our general model for the two-particle reversible $\textrm{A} + \textrm{B} \rightleftharpoons \textrm{C}$ reaction as the following
\begin{subequations}
	\begin{align}
		\PD{p}{t}(\vx,\vy,t) &= \Lp p(\vx,\vy,t) -\kp(\vx,\vy)p(\vx,\vy,t) + \int_{\Omega}\km(\vx,\vy\vert\vz)\pb(\vz,t) \, d\vz, \label{eq:FPEfwdRx} \\
		\PD{\pb}{t}(\vz,t) &= \Lpb \pb(\vz,t) -\km(\vz)\pb(\vz,t) + \int_{\Omega^2}\kp(\vz\vert\vx,\vy,t)p(\vx,\vy,t) \, d\vx \, d\vy. \label{eq:FPEbwdRx}
	\end{align}
	\label{eq:FPErevRxModel}
\end{subequations}
Here we assume a reflecting zero Neumann boundary condition on $\partial \Omega$ in each coordinate respectively (i.e. $\vx$, $\vy$ and $\vz$),
\begin{align*}
  \Dmat\paren{\nabla_{\vx,\vy}p(\vx, \vy, t)+p(\vx, \vy, t)\nabla\Phi_{\vx,\vy}} \cdot \veta(\vx,\vy) &= 0, &&(\vx,\vy) \in \partial\Omega^2,\\
   \Dbmat\paren{\nabla_{\vz} \pb(\vz,t)+\pb(\vz,t)\nabla\Phi_{\vz}} \cdot \veta_{\textrm{b}}(\vz) &= 0, &&\vz \in \partial \Omega,
\end{align*}
where $\veta(\vx,\vy)$ denotes the unit outward normal to the two-particle phase-space boundary at $(\vx,\vy)$ and $\veta_{\textrm{b}}(\vz)$ denotes the unit outward normal to $\partial\Omega$ at $\vz$. Finally, we assume the initial conditions
\begin{align*}
  p(\vx,\vy,0) &= p_0(\vx,\vy), & \pb(\vz,0) = \pbO(\vz),
\end{align*}
where $p_0$ and $\pbO$ define a proper probability distribution so that
\begin{equation*}
  \int_{\Omega^2} p_0(\vx,\vy) \, d\vx \, d\vy + \int_{\Omega} \pbO(\vz) \, d\vz = 1.
\end{equation*}
Integrating~\eqref{eq:FPErevRxModel} and using the definitions of $\km(\vz)$ and $\kp(\vx,\vy)$, this normalization of the initial conditions immediately implies that probability is conserved for all times:
\begin{equation*}
  \int_{\Omega^2} p(\vx,\vy,t) \, d\vx \, d \vy
  + \int_{\Omega} \pb(\vz,t) \, d\vz = 1.
\end{equation*}

\subsection{Detailed Balance for the Continuous Model}
\label{sec:DBContModel}
We now consider detailed balance of the reactive fluxes, and demonstrate that it allows the determination of one of the association and dissociation reactive interaction kernels in terms of the other. Assuming we are given the association kernel, $\kp(\vz | \vx, \vy)$, we derive a dissociation reaction kernel,
$\km(\vx, \vy | \vz)$, that is consistent with detailed balance.

We denote by $\peq(\vx,\vy)$ and $\pbeq(\vz)$ the steady-state solutions to \eqref{eq:FPEfwdRx} and \eqref{eq:FPEbwdRx} respectively, satisfying
\begin{equation}
	\begin{aligned}
		0 &= \Lp p(\vx,\vy,t) -\kp(\vx,\vy)p(\vx,\vy,t) + \int_{\Omega}\km(\vx,\vy\vert\vz)\pb(\vz,t) \, d\vz, \\
		0 &= \Lpb \pb(\vz,t) -\km(\vz)\pb(\vz,t) + \int_{\Omega^2}\kp(\vz\vert\vx,\vy,t)p(\vx,\vy,t) \, d\vx \, d\vy,
	\end{aligned}
	\label{eq:steadystateEqns}
\end{equation}
with reflecting zero Neumann boundary conditions on $\partial\Omega^2$ and $\partial\Omega$. As discussed in \cite{ZhangIsaacson2022}, we expect that the steady state of the system is an equilibrium state, in which the principle of detailed balance holds for the reactive fluxes, i.e.
\begin{equation}
	\kp(\vz\vert\vx,\vy)\peq(\vx,\vy) = \km(\vx,\vy\vert\vz)\pbeq(\vz).
	\label{eq:DBsteadystaterxterms}
\end{equation}
Given~\eqref{eq:DBsteadystaterxterms}, we now derive analytical expressions for the equilibrium densities, which will then allow us to derive an explicit relation defining the dissociation reaction kernel in terms of the association kernel.

Integrating \eqref{eq:DBsteadystaterxterms} with respect to $\vz$ (respectively $(\vx,\vy)$) we find
\begin{equation*}
	\begin{aligned}
		&\kp(\vx,\vy)\peq(\vx,\vy) = \int_{\Omega}\km(\vx,\vy\vert\vz)\pbeq(\vz) \, d\vz, \\
		&\km(\vz)\pbeq(\vz) = \int_{\Omega^2}\kp(\vz\vert\vx,\vy)\peq(\vx,\vy) \, d\vx \, d\vy,
	\end{aligned}
\end{equation*}
which from \eqref{eq:steadystateEqns} implies that $\Lp\peq = 0$ on $\Omega^2$ and $\Lpb\pbeq = 0$ on $\Omega$. Together with the reflecting zero Neumann boundary conditions on $\partial\Omega^2$ and $\partial\Omega$, we have
\begin{equation}
	\begin{aligned}
		\peq(\vx,\vy) &= \frac{\piAB}{\ZA \ZB}e^{-\phiA(\vx)-\phiB(\vy)}, \\
		\pbeq(\vz) &= \frac{\piC}{\ZC}e^{-\phiC(\vz)},
	\end{aligned}
	\label{eq:DBrxterms}
\end{equation}
where the normalization factors (i.e. partition functions) are given, for $S \in \{\textrm{A},\textrm{B},\textrm{C}\}$, by
\begin{equation*}
		Z_{\textrm{S}} = \int_{\Omega} e^{-\phi^{\textrm{S}}(\vx)} \, d\vx.
\end{equation*}
$\piAB$ and $\piC$ represent the probabilities of being in the unbound or bound states. Using \eqref{eq:DBrxterms} in~\eqref{eq:DBsteadystaterxterms}, we see that for the system to be consistent with detailed balance, $\kp(\vz\vert\vx,\vy)$ and $\km(\vx,\vy\vert\vz)$ must be chosen such that
\begin{equation}
	\kp(\vz\vert\vx,\vy) = \frac{\piC \ZA \ZB}{\piAB \ZC} \km(\vx,\vy\vert\vz)
	 e^{\phiA(\vx) + \phiB(\vy) - \phiC(\vz)}.
	 \label{eq:contVRkp}
\end{equation}

Let the (non-dimensional) dissociation constant $\Kd$ be the ratio of the stationary probabilities that the system is in the unbound state versus the bound state \cite{FrohnerNoe2018}, so that
\begin{equation*}
		\Kd := \frac{\piAB}{\piC} = \frac{\int_{\Omega^2}\peq(\vx,\vy) \, d\vx \, d\vy}{\int_{\Omega}\pbeq(\vz) \, d\vz}.
\end{equation*}
Given the association reactive interaction function, $\kp(\vz\vert\vx,\vy)$ in~\eqref{eq:contVRkp}, the dissociation reactive interaction function should then be defined as
\begin{equation}
	\km(\vx,\vy\vert\vz) = \Kd\frac{\ZC}{\ZA\ZB}e^{-\phiA(\vx) - \phiB(\vy) + \phiC(\vz)}\kp(\vz\vert\vx,\vy)
	\label{eq:unbindRateFcn1}
\end{equation}
to be consistent with detailed balance of the reactive fluxes holding. Note that one could equivalently define the association rate function in terms of the dissociation function, or modify both from the purely-diffusive case (as opposed to incorporating the potential terms into just one as in~\eqref{eq:unbindRateFcn1}). In the remainder we assume $\kp$ and $\km$ are chosen such that~\eqref{eq:unbindRateFcn1} holds.

\subsection{Discretization of Reaction Terms}
\label{sec:FVM_RXs}
We now summarize how we can create a jump process approximation to the reaction terms in~\eqref{eq:FPEfwdRx} and~\eqref{eq:FPEbwdRx} using the finite volume discretization developed in~\cite{IsaacsonZhang17}. We first discretize $\Omega$ into a polygonal mesh of $N$ voxels labeled by $V_i$, $i \in \{1, \dots, N\}$ with centroids $\{\vx_i\}_{i = 1\dots N}$. We define $V_{ij} = V_i \times V_j$ and $V_{ijk} = V_i \times V_j \times V_k$ with corresponding centroids labeled by $(\vx_i,\vy_j)$ and $(\vx_i,\vy_j,\vz_k)$. With these definitions, the well-mixed approximation of the probability the system is in the unbound state with the \textrm{A} molecule in $V_i$ and the \textrm{B} molecule in $V_j$ at time $t$ is given by
\begin{equation}
\label{eq:Pij}
	P_{ij}(t) = \int_{V_{ij}}p(\vx,\vy,t) \, d\vx \, d\vy \approx p(\vx_i,\vy_j,t)\abs{V_{ij}}.
\end{equation}
Similarly, the well-mixed approximation of the probability the system is in the bound state with the \textrm{C} molecule in $V_k$ at time $t$ is given by
\begin{equation}
\label{eq:Pbk}
	P_{\textrm{b}k}(t) = \int_{V_k}\pb(\vz,t) \, d\vz \approx \pb(\vz_k,t)\abs{V_k}.
\end{equation}

We will subsequently drop the spatial transport terms in \eqref{eq:FPEfwdRx} and \eqref{eq:FPEbwdRx} as we approximate them through the EAFE method of Section~\ref{sec:DriftDiffApprox}. To apply the finite volume discretization to the reaction terms as in~\cite{IsaacsonZhang17}, we integrate both sides of \eqref{eq:FPEfwdRx} and \eqref{eq:FPEbwdRx} over $V_{ij}$ and $V_{k}$ respectively. \eqref{eq:FPEfwdRx} is then approximated by
\begin{equation}
	\begin{aligned}
		\D{P_{ij}}{t} &\approx -\frac{1}{\abs{V_{ij}}}P_{ij}(t)\int_{V_{ij}}\kp(\vx,\vy) \, d\vx \, d\vy + \sum_k\frac{1}{\abs{V_k}}P_{\textrm{b}k}(t)\int_{V_{ijk}}\km(\vx,\vy\vert\vz) \, d\vx \, d\vy \, d\vz \\
		&= -\kp_{ij}P_{ij}(t) + \sum_k\km_{ijk}P_{\textrm{b}k}(t),
	\end{aligned}
	\label{eq:FPEfwdRxFV}
\end{equation}
where
\begin{subequations}
\begin{align}
	\kp_{ij} &:= \frac{1}{\abs{V_{ij}}}\int_{V_{ij}}\kp(\vx,\vy) \, d\vx \, d\vy, \label{eq:kpij}\\
	\km_{ijk} &:= \frac{1}{\abs{V_k}}\int_{V_{ijk}}\km(\vx,\vy\vert\vz) \, d\vx \, d\vy \, d\vz. \label{eq:kmijk}
\end{align}
\end{subequations}
$\kp_{ij}$ defines the probability per unit time that an \textrm{A} molecule in $V_i$ and a \textrm{B} molecule in $V_j$ react to produce a \textrm{C} molecule in $\Omega$. Likewise, $\km_{ijk}$ can be interpreted as the probability per unit time that a \textrm{C} molecule in $V_k$ unbinds to create an \textrm{A} molecule in $V_i$ and a \textrm{B} molecule in $V_j$.

By the same approach, the reaction terms in \eqref{eq:FPEbwdRx} are approximated by
\begin{equation}
	\begin{aligned}
		\D{P_{\textrm{b}k}}{t} &\approx -\frac{1}{\abs{V_k}}P_{\textrm{b}k}(t)\int_{V_k}\km(\vz) \, d\vz + \sum_{i,j}\frac{1}{\abs{V_{ij}}}P_{ij}(t)\int_{V_{ijk}}\kp(\vz\vert\vx,\vy) \, d\vx \, d\vy \, d\vz \\
		&= -\km_{k}P_{\textrm{b}k}(t) + \sum_{i,j}\kp_{ijk}P_{ij}(t),
	\end{aligned}
	\label{eq:FPEbwdRxFV}
\end{equation}
where
\begin{subequations}
\begin{align}
	\km_{k} &:= -\frac{1}{\abs{V_k}}\int_{V_k}\km(\vz) \, d\vz, \label{eq:kmk}\\
	\kp_{ijk} &:= \frac{1}{\abs{V_{ij}}}\int_{V_{ijk}}\kp(\vz\vert\vx,\vy) \, d\vx \, d\vy \, d\vz. \label{eq:kpijk}
\end{align}
\end{subequations}
$\km_{k}$ gives the probability per unit time that a \textrm{C} molecule in $V_k$ unbinds into \textrm{A} and \textrm{B} molecules within $\Omega$. Similarly, $\kp_{ijk}$ defines the probability per unit time that an \textrm{A} molecule in $V_i$ and a \textrm{B} molecule in $V_j$ react and create a \textrm{C} molecule in $V_k$.

Using the definitions of $\kp(\vx,\vy)$ and $\km(\vz)$, we can rewrite $\kpij$ and $\kmk$ as
\begin{subequations}
	\begin{align}
		\kpij &= \sum_{k}\kpijk, \\
		\kmk &= \sum_{i,j}\kmijk.
	\end{align}
\end{subequations}

\subsection{Unstructured Mesh CRDDME with Reversible Reactions}
\label{sec:discreteModel}
To arrive at the final unstructured mesh CRDDME with reversible $\textrm{A} + \textrm{B} \rightleftharpoons \textrm{C}$ reaction, we combine the EAFE discretization~\eqref{eq:semidisFEMProb} and~\eqref{eq:discreteDDOperator} for the spatial transport operators in~\eqref{eq:FPEfwdRx} and~\eqref{eq:FPEbwdRx} with the finite volume discretizations~\eqref{eq:FPEfwdRxFV} and~\eqref{eq:FPEbwdRxFV} for the reaction operators \eqref{eq:kmijk} and \eqref{eq:kpijk} summarized in the last section. Let $P_{ij}(t)$ denote the probability for the \textrm{A} and \textrm{B} molecules to be in $V_{ij}$ at time $t$, and $P_{\textrm{b}k}(t)$ the probability for the \textrm{C} molecule to be in $V_k$ at time $t$. We arrive at the semi-discrete jump process CRDDME approximation
\begin{subequations} \label{eq:twopartCRDDME}
  \begin{equation}
    \label{eq:twopartCRDDMEfwd}
    \begin{aligned}
      \D{P_{ij}}{t} &= \sum_{i'=1}^{N} \brac{(\diffophA)_{i i'} P_{i'j}(t)
        - (\diffophA)_{i' i} P_{ij}(t)} \\
      &+ \sum_{j'=1}^{N} \brac{ (\diffophB)_{j j'} P_{ij'}(t) - (\diffophB)_{j' j} P_{ij}(t)}
      -\sum_{k = 1}^N\kpijk P_{ij}(t) + \sum_{k=1}^{N} \kmijk P_{\textrm{b}k}(t),
    \end{aligned}
  \end{equation}
  \begin{multline}
    \label{eq:twopartCRDDMEbwd}
    \D{P_{\textrm{b}k}}{t} = \sum_{k'=1}^{N} \brac{(\diffophC)_{k k'} P_{\textrm{b}k'}(t)
      - (\diffophC)_{k' k} P_{\textrm{b}k}(t)}
    -\sum_{i,j = 1}^N\kmijk P_{\textrm{b}k}(t)  + \sum_{i,j=1}^{N} \kpijk P_{ij}(t),
  \end{multline}
\end{subequations}
where $\diffophA$, $\diffophB$, and $\diffophC$ denote the EAFE discretization of each particle's associated drift-diffusion operator (with potentials $\phiA$, $\phiB$ and $\phiC$ respectively) given by ~\eqref{eq:semidisFEMProb} and~\eqref{eq:discreteDDOperator}.

The jump process reaction approximation can be simulated in several equivalent ways. One approach uses the transition rates $\kpij$ and $\kmk$ to determine which reaction occurs next, and then sample the reaction placement probabilities
\begin{subequations}
	\begin{align}
		\kpkij &:= \frac{\kpijk}{\kpij} \\
		\kmkij &:= \frac{\kmijk}{\kmk}
	\end{align}
\end{subequations}
to determine where products should be placed. Alternatively, one can directly use the transition rates $\kpijk$ and $\kmijk$ to simultaneously sample both when reactions occur and where products are placed. Note that $\kpijk$ and $\kmijk$ can be pre-tabulated by numerical integration prior to running stochastic simulation algorithm (SSA) simulations. If $\kp(\vx,\vy)$ is sufficiently smooth, standard Gaussian quadrature approximations can be used~\cite{IsaacsonZhang17}. When $\kp(\vx,\vy)$ is given by a discontinuous Doi interaction, hyper-volume intersections can be calculated exactly to both reduce the dimensionality of the needed integrals and ensure a continuous integrand in any quadrature approximations~\cite{IsaacsonCRDME2013}.

\begin{table}[tb]
	\begin{tabular}{|c|c|c|c|}
	\hline
	\multicolumn{1}{|l|}{}              & Transitions                                                     & Transition Rates          & Upon Transition Event                                             \\ \hline
	\multirow{3}{*}{Spatial Movement}   & $\textrm{A}_j \to \textrm{A}_i$                                 & $(\diffophA)_{ij}a_j$     & $\textrm{A}_i := \textrm{A}_i+1$, $\textrm{A}_j:=\textrm{A}_j-1$, \\ \cline{2-4}
										& $\textrm{B}_j \to \textrm{B}_i$                                 & $(\diffophB)_{ij}b_j$     & $\textrm{B}_i := \textrm{B}_i+1$, $\textrm{B}_j:=\textrm{B}_j-1$, \\ \cline{2-4}
										& $\textrm{C}_j \to \textrm{C}_i$                                 & $(\diffophC)_{ij}c_j$     & $\textrm{C}_i := \textrm{C}_i+1$, $\textrm{C}_j:=\textrm{C}_j-1$, \\ \hline
	\multirow{4}{*}{Chemical Reactions} & \multirow{2}{*}{$\textrm{A}_i + \textrm{B}_j \to \textrm{C}_k$} & \multirow{2}{*}{$\kpijk a_i b_j$} & $\textrm{A}_i := \textrm{A}_i-1$, $\textrm{B}_j:=\textrm{B}_j-1$, \\
										&                                                                 &                           & $\textrm{C}_k := \textrm{C}_k+1$.                                 \\ \cline{2-4}
										& \multirow{2}{*}{$\textrm{C}_k \to \textrm{A}_i + \textrm{B}_j$} & \multirow{2}{*}{$\kmijk c_k$} & $\textrm{C}_k := \textrm{C}_k-1$,                                 \\
										&                                                                 &                           & $\textrm{A}_i := \textrm{A}_i+1$, $\textrm{B}_j:=\textrm{B}_j+1$, \\ \hline
	\end{tabular}
	\caption{Summary of spatial and chemical reactions for the jump process approximation of the multiparticle $\textrm{A} + \textrm{B} \leftrightarrows \textrm{C}$ reaction system. Here $a_i$ denotes the number of $\textrm{A}$ particles in voxel $V_i$, with $b_j$ and $c_k$ defined similarly. Transition rates give the probability per time for a transition to occur (i.e. propensities). The final column explains how to update the system upon occurrence of a transition event. The transition rate matrix $\mathcal{L}_h$ and reaction rates $\kappa^{\pm}_{i j k}$ are given by \eqref{eq:discreteDDOperator}, \eqref{eq:kpijk}, and \eqref{eq:kmijk} respectively.}
	\label{Tab:CRDDMETransitionRates}
\end{table}	
While \eqref{eq:twopartCRDDME} describes a system containing only one $\textrm{A}$ and one $\textrm{B}$ molecule, or one $\textrm{C}$ molecule in the bound state, it is straightforward to generalize the system to include arbitrary numbers of molecules of each species. In Appendix~\ref{app:multipartModel} we formulate the corresponding continuous model for such multiparticle systems by generalizing the two-particle system given by \eqref{eq:FPErevRxModel}. Using the same approach as described in \cite{IsaacsonZhang17}, we derive a general CRDDME for multiparticle systems where particles of the same species experience the same one-body potential \eqref{S1eq:masterEqnProb}. In Table~\ref{Tab:CRDDMETransitionRates} we summarize the resulting set of spatial and chemical transitions, along with associated transition rates. Here we note that the transition rates for the two-particle system fully determine the transition rates for the multi-particle system. The transitions given by Table~\ref{Tab:CRDDMETransitionRates} can be interpreted as a large collection of chemical reactions, and as such exact realizations of the jump processes associated with the CRDDME for general multiparticle reaction systems can be simulated using any SSA~\cite{MarchettiSimCompBio}. In the remainder, for all simulations that include reactions we use the Next Reaction Method SSA of Gibson and Bruck~\cite{GibsonBruckJPCHEM2002} to generate exact realizations of the jump processes associated with the CRDDME. All reported statistics for problems involving reactions are then sample statistics, estimated from running many such simulations.

\subsection{Detailed Balance in the CRDDME}
\label{sec:DBCRDDME}
In the case of reversible reactions, such as the $\textrm{A} + \textrm{B} \leftrightarrows \textrm{C}$ reaction we consider here, we would like the CRDDME to also be consistent with detailed balance of the reaction fluxes holding at equilibrium. This means that we would like a discrete version of the spatially-continuous equilibrium state~\eqref{eq:DBrxterms}, with steady-state probabilities $\equil{P}_{i j}$ and $\equil{P}_{\textrm{b}k}$, to be a steady-state of the CRDDME which satisfies the detailed balance condition that
\begin{equation}
	\kpijk\equil{P}_{ij} = \kmijk\equil{P}_{\textrm{b}k}.
	\label{eq:discreteDB}
\end{equation}
When this condition holds, we find that
\begin{align*}
	0 &= \sum_{i'=1}^{N} \brac{(\diffophA)_{i i'} \equil{P}_{i'j} - (\diffophA)_{i'i} \equil{P}_{ij}}
     + \sum_{j'=1}^{N} \brac{(\diffophB)_{j j'} \equil{P}_{ij'} - (\diffophB)_{j'j} \equil{P}_{ij}} \\
	0 &= \sum_{k'=1}^{N} \brac{(\diffophC)_{k k'} \equil{P}_{\textrm{b}k'}
	 - (\diffophC)_{k' k} \equil{P}_{\textrm{b}k}}
\end{align*}
which in turn imply that
\begin{align*}
	\equil{P}_{i j} &= \frac{\pihatAB}{\ZhatA \ZhatB} e^{-\phiA_i - \phiB_j} \abs{V_{i j}}, \\
	\equil{P}_{\textrm{b} k} &= \frac{\pihatC}{\ZhatC} e^{-\phiC_k} \abs{V_k},
\end{align*}
where $\hat{Z}^{\textrm{S}}$ for $\textrm{S} \in \{\textrm{A}, \textrm{B}, \textrm{C}\}$ is the mesh partition function,
\begin{equation*}
	\hat{Z}_{\textrm{S}} = \sum_j e^{-\phi^{\textrm{S}}_j} \abs{V_j},
\end{equation*}
and $\pihatAB$ and $\pihatC$ denote the probabilities of being in the bound or unbound states. Substituting $\equil{P}_{i j}$ and $\equil{P}_k$ into~\eqref{eq:discreteDB}, we can define one of $\kpijk$ or $\kmijk$ in terms of the other, i.e.
\begin{equation*}
	\kmijk = \frac{\pihatAB \ZhatC}{\pihatC \ZhatA \ZhatB} \frac{\abs{V_{i j}}}{\abs{V_k}}
		e^{\phiC_k - \phiA_i - \phiB_j} \kpijk.
\end{equation*}
As in the spatially continuous case, we assume the ratio of the probabilities of being in the unbound to bound state is given by the dissociation constant of the reaction, i.e. that
\begin{equation} \label{eq:discDBkmfromkp}
	\kmijk = \frac{\Kd \ZhatC}{\ZhatA \ZhatB} \frac{\abs{V_{i j}}}{\abs{V_k}}
		e^{\phiC_k - \phiA_i - \phiB_j} \kpijk.
\end{equation}
Then, in the case that the potentials are constant, so that particles move purely by diffusion, the value of $\kmijk$ obtained from~\eqref{eq:discDBkmfromkp} is identical to the formula~\eqref{eq:kmijk} when the continuous interaction functions of the volume reactivity model are chosen to satisfy the detailed balance relation~\eqref{eq:unbindRateFcn1}~\cite{IsaacsonZhang17}.

Unfortunately, in the case of non-constant potentials, \eqref{eq:discDBkmfromkp} and the formula~\eqref{eq:kmijk} for $\kmijk$ obtained by directly discretizing the continuous unbinding kernel are no longer identical. That is, applying the finite volume discretization to both the reaction terms does not lead to reactive transitions rates that are consistent with detailed balance at the discrete level. To ensure our discrete model satisfies detailed balance, we choose $\kmijk$ using~\eqref{eq:discDBkmfromkp} assuming $\kpijk$ is given by the discretized formula~\eqref{eq:kpijk}. We note that this is equivalent to applying a quadrature rule to \eqref{eq:kmijk} under the assumption that the continuous kernels satisfy the detailed balance relation~\eqref{eq:unbindRateFcn1}. That is
\begin{align*}
	\kmijk &= \frac{1}{\abs{V_k}} \int_{V_{i j k}} \km(\vx,\vy\vert\vz) \, d\vx \, d\vy \, d\vz \\
	&= \frac{\Kd \ZC}{\ZA \ZB} \frac{1}{\abs{V_k}} \int_{V_{i j k}}
	e^{-\phiA(\vx) - \phiB(\vy) + \phiC(\vz)} \kp(\vz \vert \vx,\vy) \, d\vx \, d\vy \, d\vz \\
	&\approx \frac{\Kd \ZhatC}{\ZhatA \ZhatB} \frac{1}{\abs{V_k}}
	e^{\phiC_k - \phiA_i - \phiB_j} \int_{V_{i j k}} \kp(\vz \vert \vx, \vy) \, d\vx \, d\vy \, d\vy \\
	&= \frac{\Kd \ZhatC}{\ZhatA \ZhatB}  \frac{\abs{V_{i j}}}{\abs{V_k}}
	e^{\phiC_k - \phiA_i - \phiB_j} \kpijk.
\end{align*}

We therefore have a choice. We can use the $\kmijk$ that is obtained from direct finite volume approximation of the spatially-continuous model, assuming detailed balance of the continuous interaction functions, and obtain a formula that does not satisfy the discrete detailed balance condition. Alternatively, we can use the discrete detailed balance condition~\eqref{eq:discDBkmfromkp} to define one of $\kpijk$ and $\kmijk$ in terms of the other, thereby making a second implicit quadrature-type approximation, but ensuring consistency with the discrete detailed balance condition. In the remainder we use this latter approach, demonstrating in the numerical examples that the statistics we examine demonstrate an empirical second order rate of convergence with this additional approximation, consistent with what might be conjectured should the underlying CRDDME solution be converging at second order in $L^2$ to the solution of the PIDEs of the VR model~\cite{HeldmanTheis2023}.

The two-particle detailed balance condition given by \eqref{eq:discreteDB} can be generalized to the multiparticle system. In Appendix ~\ref{app:multipartDB} we demonstrate the consistency of detailed balance at the discrete level for the multiparticle system. We derive the corresponding detailed balance condition for the discrete multiparticle system and subsequently show that the multiparticle detailed balance condition can be satisfied whenever the two-particle detailed balance condition \eqref{eq:discreteDB} is enforced.

\section{Numerical Examples}
\label{sec:numerical_examples}
We now illustrate the convergence and accuracy of the unstructured mesh CRDDME with reactions through several examples. For all simulations we generate exact realizations of the jump process defined in Section \ref{sec:discreteModel}, using the next reaction method SSA \cite{GibsonBruckJPCHEM2002}. In Section \ref{sec:annihilRx} we demonstrate that several reaction time statistics converge to finite values as the mesh size approaches zero for the two-particle $\textrm{A} + \textrm{B} \to \emptyset$ annihilation reaction within a circle. We examine the standard discontinuous Doi interaction. With convergence established for the forward reaction approximation, we then confirm in Section \ref{sec:revRx} that statistics of the two-particle reversible $\textrm{A} + \textrm{B} \rightleftharpoons \textrm{C}$ reaction converge to the solution of the Doi model by comparing with Brownian Dynamics (BD) simulations.

These numerical examples demonstrate (empirical) convergence of statistics associated with our hybrid discretization approach that combines a finite element approximation for spatial transport and a finite volume approximation for reactive components. Discussion of related work on hybrid PDE discretizations, rigorous convergence results for hybrid finite element / finite volume CRDME-type PDE discretizations, and rigorous results for convergence of such CRDME-type models to the underlying VR model are available in~\cite{HeldmanTheis2023,Heldman2024}.

\subsection{$\textrm{A} + \textrm{B} \to \emptyset$ Annihilation Reaction}
\label{sec:annihilRx}
We begin by examining the annihilation reaction, $\textrm{A} + \textrm{B} \to \emptyset$ in a system with just one A molecule and one B molecule in two-dimensions ($d = 2$) with a Doi reaction mechanism
\begin{equation}
	\kp(\vx,\vy) = \lambda\ind_{B_{\rb}(\vO)}(\abs{\vx - \vy}),
\end{equation}
where two molecules attempt to react with probability per unit time $\lambda$ when within $\varepsilon$.
We assume each molecule diffuses either within a square centered at the origin with side length $0.1$ $\mu$m (i.e. $[-0.05, 0.05]\times[-0.05, 0.05]$) or a circle centered at $(0.05,0.05)$ with radius $0.1$ $\mu$m. A diffusion constant of $10$ $\mu$m$^{2}$s$^{-1}$ is used for both \textrm{A} and \textrm{B} molecules. Each molecule is assumed to experience the same background potential, $\phi^A(\vx) = \phi^B(\vy) = \phi(\vv)$, given by
\begin{equation}
	\phi(\vv) = x^2 + y^2, \,\,\,\,\,\, \vv = (x,y).
	\label{eq:ExPotential}
\end{equation}
In the absence of reaction, this background potential drives particles to the origin, where $\phi(\vv)$ admits its minimum value. The reaction radius, $\rb$, is chosen to be 1 nm, and the reaction rate, $\lambda$, is set to be $10^9$ s$^{-1}$.

\begin{figure}[!tbp]
  \centering
  \includegraphics[width = 1\textwidth]{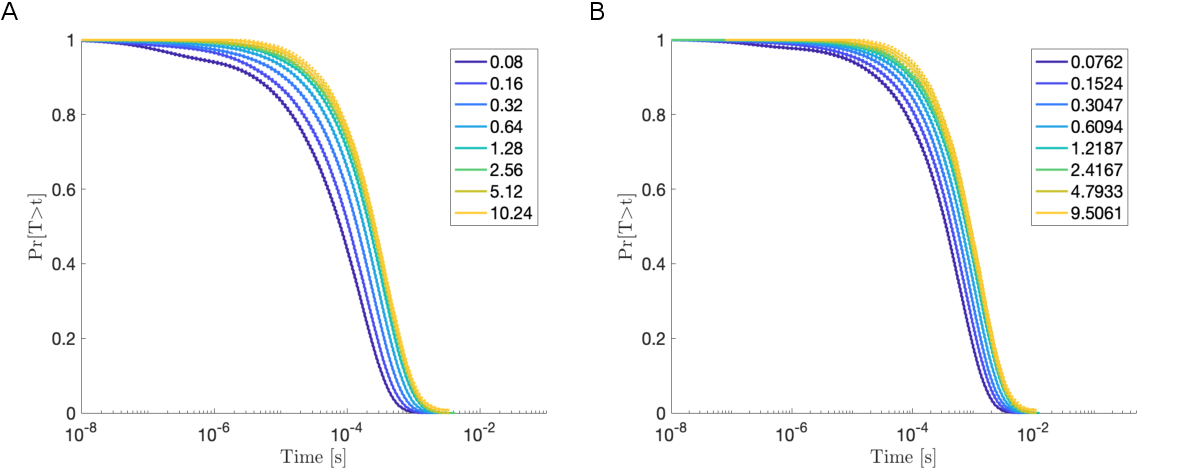}
  \caption{Survival time distributions vs. $t$ from the two-particle
  $\textrm{A} + \textrm{B} \to \varnothing$ reaction in a square in panel A and in a disk in panel B. In both cases each curve was estimated from 100000 simulations using the \texttt{ecdf} command in MATLAB. The legends give the ratio, $\rb$/$h$ as the mesh is refined ($h$ is approximately successively halved). See Section~\ref{sec:annihilRx} for parameter values.}
  \label{fig:sq_cr_SurvProb}
\end{figure}

Let $\Tb$ denote the time for the two molecules to react when each starts uniformly distributed in $\Omega$. The corresponding survival time distribution, $\Pr[\Tb > t]$, is given by
\begin{equation*}
	\Pr[\Tb > t] = \int_{\Omega} p(\vx,\vy,t) \, d\vx \, d\vy,
\end{equation*}
where $p(\vx,\vy,t)$ satisfies \eqref{eq:FPEfwdRx} with $\km(\vx,\vy\vert\vz) = 0$. We estimate the survival time distribution from the numerically sampled reaction times using the \texttt{ecdf} command in MATLAB. In Fig.~\ref{fig:sq_cr_SurvProb} we demonstrates the convergence (to within sampling error) of the estimated survival time distribution of the unstructured mesh CRDDME using a Doi interaction. The survival time distributions are seen to converge as the maximum mesh width $h \to 0$.

\begin{figure}[!tbp]
  \centering
  \includegraphics[width = 1\textwidth]{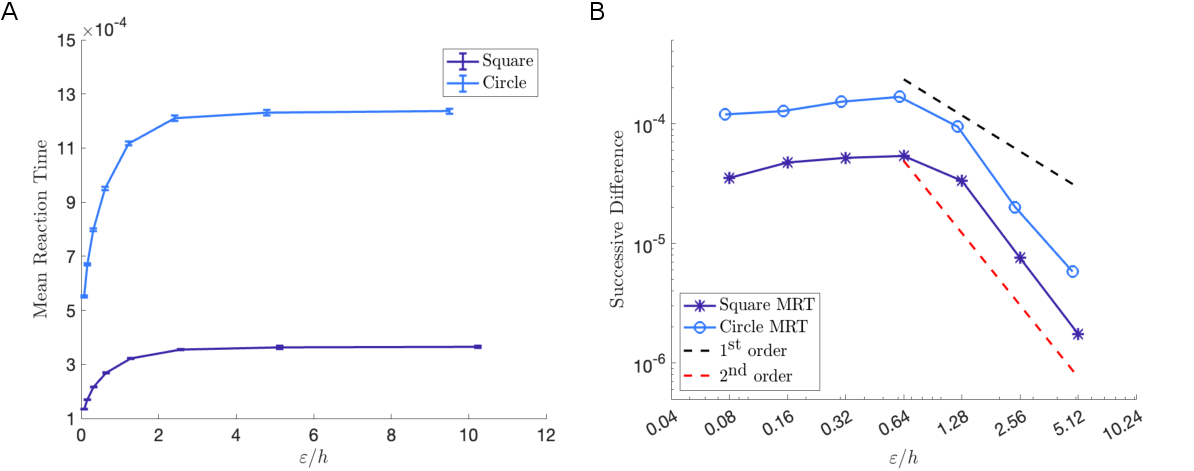}
  \caption{Convergence of the mean reaction time $\avg{\Tb}$. In panel A we plot the mean reaction time $\avg{\Tb}$ vs. $\rb$/$h$ as $h$ is (approximately) successively halved. Each mean reaction time for the two cases was estimated from 100000 simulations. 95$\%$ confidence intervals are drawn on each data point, but for some points are smaller than the marker labeling the point. See Section~\ref{sec:annihilRx} for parameter values. In panel B we demonstrate the rate of convergence when using a square and a disk domain by plotting the difference between successive points on the $\avg{\Tb}$ vs $\rb/h$ curves from panel A. The smaller of the two $h$ values is used for labeling. The effective convergence rate to zero of the FPE with both reaction mechanisms scales roughly like $O(h^2)$. Note that the search time in the square is smaller in part as the area of the circle is a factor of $\pi$ larger than the square.}
    \label{fig:convergMRT}
\end{figure}

To study the rate of convergence we also examined the mean reaction time $\E[\Tb]$, defined by
\begin{equation*}
	\E[\Tb] = \int_0^\infty \Pr[\Tb > t] \, dt.
\end{equation*}
We estimated the mean reaction time from the numerically sampled reaction times by calculating the sample mean. In Fig.~\ref{fig:convergMRT}A we show the sample mean reaction times for the Doi reaction as $\rb/h$ is varied. We see that as $\rb/h$ is increased (i.e. $h\to0$) the sample mean reaction times converge to a finite value. Figure~\ref{fig:convergMRT}B illustrates the rate of convergence using the Doi interaction by plotting the successive difference of the estimated mean reaction times as $h$ is decreased (approximately halved). For $h$ sufficiently small, the empirical rate of convergence for both reaction mechanisms is roughly second order.

\subsection{\textrm{A} + \textrm{B} $\rightleftharpoons$ \textrm{C} Reversible Reaction}
\label{sec:revRx}
We now consider the bimolecular reversible $\textrm{A} + \textrm{B} \rightleftharpoons \textrm{C}$ reaction in a system initialized with only one \textrm{C} molecule. The corresponding volume reactivity drift-diffusion model is given by \eqref{eq:FPEfwdRx} and \eqref{eq:FPEbwdRx}. All molecules are assumed to have the same diffusion constant, $\DA = \DB = \DC = 0.1 \mu$m$^2$s$^{-1}$, and experience the same background potential $\phiA(\vx) = \phiB(\vy) = \phiC(\vz) = \phi(\vv)$ defined by \eqref{eq:ExPotential}. The domain $\Omega$ is chosen to be a circle centered at $(0.05,0.05)$ with radius $0.1$ $\mu$m. We assume a reflecting Neumann boundary condition on $\partial\Omega$ in each of $\vx$, $\vy$, and $\vz$ coordinates. For the forward $\textrm{A} + \textrm{B} \to \textrm{C}$ we use the Doi reaction model
\begin{equation}
	\kp(\vz\vert\vx, \vy) = \lambda\ind_{B_{\rb}(\vO)}(\abs{\vx - \vy})\delta(\vz - \gamma\vx - (1-\gamma)\vy),
	\label{eq:fwdRxModel}
\end{equation}
with a reaction radius $\rb = 1$ nm and $\gamma = \frac{1}{2}$. Note that \eqref{eq:fwdRxModel} models the reaction of \textrm{A} and \textrm{B} particles by assuming that two particles may react with probability per time $\lambda$ when within $\varepsilon$ of each other, described by the indicator function term, and that a newly created $\textrm{C}$ particle is placed on the line connecting the \textrm{A} and \textrm{B} particles, described by the $\delta$ function term. There are several different choices for $\kp(\vz\vert\vx,\vy)$, here we choose the association kernel such that it is the same as in the pure diffusive case. The corresponding dissociation rate can be obtained using \eqref{eq:unbindRateFcn1}
 \begin{equation}
	\km(\vx,\vy\vert\vz) = \Kd\lambda\frac{Z_C}{Z_AZ_B}e^{-\phi(\vx)-\phi(\vy)+\phi(\vz)}\ind_{B_{\rb}(\vO)}(\abs{\vx - \vy})\delta(\vz - \gamma\vx - (1-\gamma)\vy).
	\label{eq:bwdRxModel}
\end{equation}
In the simulations, the corresponding $\kmijk$ is obtained by evaluating \eqref{eq:bwdRxModel} using numerical quadrature.

For all simulations the \textrm{C} molecule is initially placed randomly within $\Omega$, corresponding to the initial conditions
\begin{equation*}
	p(\vx,\vy, 0) = 0, \qquad \pb(\vz,0) = \frac{1}{\abs{\Omega}}.
\end{equation*}

\begin{figure}[!tbp]
  \centering
  \includegraphics[width=\textwidth]{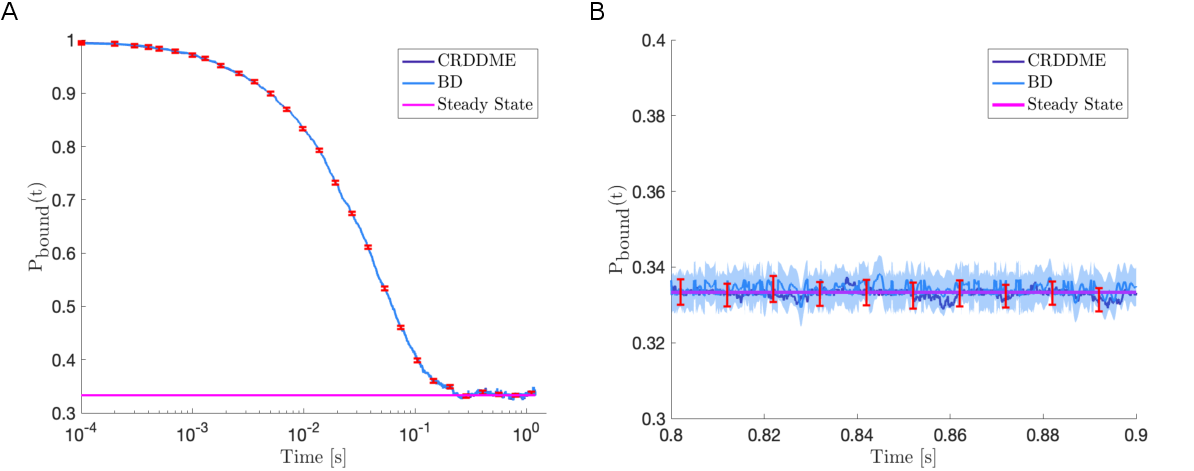}
  \caption{Probability that molecules are in bound state $\Pbound(t)$ versus time for CRDDME SSA simulations and BD simulations. Dark purple curves correspond to $\Pbound(t)$ from CRDDME simulations and the light blue to $\Pbound(t)$ from BD simulations. The magenta solid line shows the steady-state $\Pbound$ value for $\Kd = 2$. $95\%$ confidence intervals for the CRDDME simulations are drawn with red error bars, and as a translucent blue ribbon for the BD simulations. Panel B shows a zoomed in version for a portion of panel A when the system reaches steady states, making clear the scale of the 95\% confidence intervals. Each curve was estimated from 100000 simulations by dividing the number of simulations in which the system is in the bound state at a given time by the total number of simulations. The domain $\Omega$ is a circle centered at $(0.05, 0.05)$ with radius $r = 0.1\mu$m and was discretized into $44945$ polygons.}
    \label{fig:convergPbBD}
\end{figure}

\begin{figure}[!tbp]
  \centering
  \includegraphics[width=\textwidth]{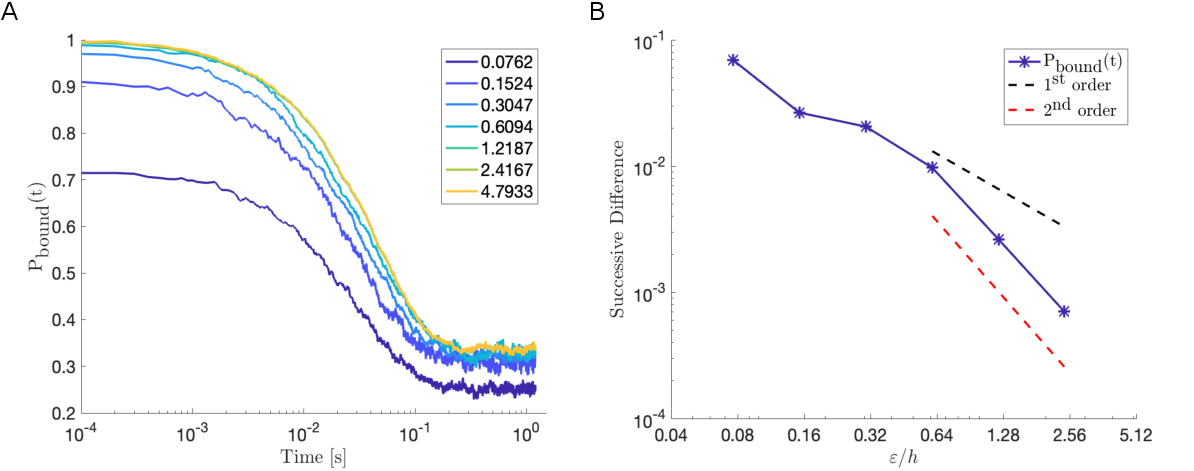}
  \caption{Convergence of the probability the molecules are in the bound state, $\Pbound(t)$, as $h \to 0$. In panel A we plot $\Pbound(t)$ vs. $\rb$/$h$ as $h$ is (approximately) successively halved. Each curve was estimated from 100000 simulations. We can observe convergence as the mesh width $h$ approaches $0$. See Section~\ref{sec:revRx} for parameter values. In panel B we demonstrate the rate of convergence by plotting the difference between successive points on the $\Pbound(t)$ vs $\rb/h$ curves from panel A at $t = 0.1$s. The smaller of the two $h$ values is used for labeling. The effective convergence rate to zero of the CRDDME with both reaction mechanisms scales roughly like $O(h^2)$.}
  \label{fig:convergPb}
\end{figure}

To confirm that the unstructured mesh CRDDME with reversible reaction converges to the solution of the drift-diffusion Doi VR model, we compare reaction statistics from SSA simulations against reaction statistics calculated from Brownian Dynamics (BD) simulations. Our BD simulations use the methods from~\cite{ReaddyPLOS2013,ReaddyPLOS2019}, modifying the diffusion update step to simulate drift-diffusion with a fixed time step of $dt = 10^{-10}$ s. For all simulations the association rate constant in \eqref{eq:fwdRxModel} is set to be $\lambda = 10^{6}$ s$^{-1}$, and the ratio of the steady-state probabilities is chosen to be $\Kd = 2.0$. We denote by $\Pbound(t)$ the probability the system is in the \textrm{C} state (i.e. the \textrm{A} and \textrm{B} molecules are bound together) at time $t$,
\begin{equation*}
	\Pbound(t) = \int_{\Omega}\pb(\vz,t) \, d\vz.
\end{equation*}
With the choice of $\Kd = 2$, in the limit as $t \to \infty$
\begin{equation*}
	\Pbound(t) = \frac{1}{1+\Kd} \approx 0.3333.
\end{equation*}
We estimated $\Pbound(t)$ numerically by averaging the number of \textrm{C} molecules in the system at a fixed time over the total number of SSA (resp. BD) simulations. We demonstrate in Fig.\ref{fig:convergPbBD} that $\Pbound(t)$ from SSA simulations agrees with $\Pbound(t)$ from BD simulations to statistical error. Figure \ref{fig:convergPb}B shows the convergence of $\Pbound(t)$ from the CRDDME simulations as the mesh width $h$ goes to zero. To demonstrate the rate of convergence, we plot the successive difference of the estimated $\Pbound(t)$ at $t = 0.1$s as $h$ is approximately halved. As $\rb/h$ is increased, the rate of convergence is approximately second order.

We note that with the choice of the potential function \eqref{eq:ExPotential}, molecules will aggregate at the origin at steady state. As such, whenever the \textrm{C} molecule unbinds, the newborn \textrm{A} and \textrm{B} will quickly bind again, which gives rise to the noisiness of the steady-state $\Pbound$. 

\subsection{Multiparticle Systems}
\label{sec:multiParticleSysm}
\begin{figure}[!tbp]
  \centering
  \includegraphics[width=\textwidth]{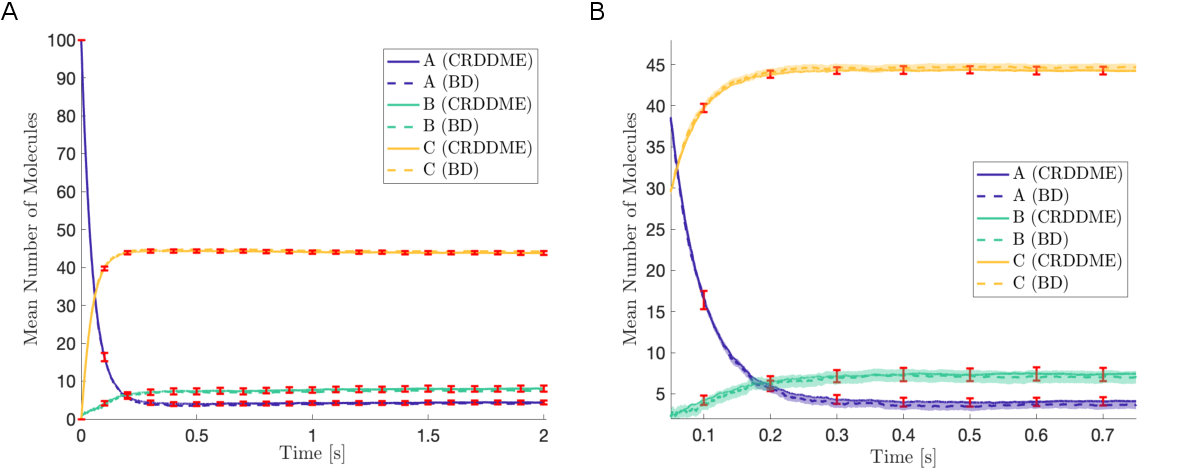}
  \caption{Mean number of molecules of each species vs. time for CRDDME SSA and BD simulations when $\rb = 1$nm. In panel A, we plot the average number of molecules that are in each state versus time for CRDDME SSA simulations and BD simulations The confidence intervals of the BD simulations are omitted here for a better visualization purpose. In panel B, we provide a zoomed in region of the figure in panel A with the corresponding $95\%$ confidence intervals of the BD simulations plotted using transparent ribbons. In both panels, solid lines correspond to CRDDME simulations whereas dashed lines correspond to BD simulations. $95\%$ confidence intervals of the CRDDME simulations are drawn in red. Each curve was estimated from 100 simulations by calculating the sample average for the number of molecules of each species at a given time. All BD simulations used a time-step of $dt = 10^{-10}$s.}
    \label{fig:multPartSysmComparison}
\end{figure}

In the previous subsections, we have demonstrated the convergence of the method for basic bimolecular reactions. We now demonstrate that the method can also accurately resolve multiparticle systems by considering the following reaction system
\begin{equation}
	\begin{aligned}
		&\textrm{A} \xrightleftharpoons[\koff]{\kon} \textrm{B} \\
  		\textrm{A} &+ \textrm{B} \xrightleftharpoons[\mu]{\lambda} \textrm{C}
	\end{aligned}
	\label{eq:MultPartSysm}
\end{equation}
where the first-order reaction rates $\kon = 10$s$^{-1}$ and $\koff = 5$s$^{-1}$, and the forward $\textrm{A} + \textrm{B} \to \textrm{C}$ association rate is set to be $\lambda = 10^5$s$^{-1}$. The $\textrm{C} \to \textrm{A} + \textrm{B}$ dissociation rate is determined implicitly by setting the ratio of the steady-state probabilities $\Kd = 2.0$. The domain $\Omega$ is chosen to be a circle of radius $0.1 \mu$m centered at $(0.5, 0.5)$, discretized into $44945$ elements. All molecules are assumed to diffuse at the same rate $\DA = \DB = \DC = 0.1 \mu$m$^2$s$^{-1}$, and the diffusive motions are affected by the same background potential $\phiA(\vx) = \phiB(\vy) = \phiC(\vz) = \phi(\vv)$ defined by \eqref{eq:ExPotential}. In Fig.~\ref{fig:multPartSysmComparison} we plot the time evolution of the average number of molecules of each species from CRDDME SSA and BD simulations. The estimated average number of molecules between the CRDDME SSA and BD when averaged over $100$ simulations agree to statistical error.
\subsection{TCR-pMHC Transport in Immune Synapse Formation}
\label{sec:ISmodel}
When T cells and antigen presenting cells (APCs) interact, binding interactions between proteins on the two cells' surfaces can result in adhesion and the formation of an immunological synapse (IS)~\cite{DharanFarago2017,Siokisetal2018}. A striking component of ISs is the segregation of several of these crosslinked proteins, with a circular core of TCR-pMHC complexes surrounded by a peripheral region of LFA1-ICAM1 complexes. An open question that a number of recent in silico studies have investigated is to understand what drives this segregation of the protein complexes~\cite{DharanFarago2017,Siokis2017,Siokisetal2018}. These studies have relied on lattice models that provide more coarse-grained representations for the underlying drift-diffusion dynamics of proteins than the VR model the CRDDME approximates. We now demonstrate that, with the exception of two-body interaction potentials, the CRDDME can handle the key physical components that form the models presented in these works. In particular, we show that its predications for the dynamics of TCR-pMHC aggregation in the central core of the IS are consistent with what is reported via in silico knock-down experiments that removed attractive two-body interactions~\cite{DharanFarago2017}.

The model we use involves T-cell receptors (TCRs) and pMHCs that react
through the reversible bimolecular reaction
\begin{equation}
	\textrm{TCR} + \textrm{pMHC} \xrightleftharpoons[\mu]{\lambda} \textrm{TCR-pMHC}.
	\label{eq:ISmodel}
\end{equation}
As the CRDDME does not yet support two-body interactions, we do not take into consideration LFA-1-ICAM-1, which primarily interacted with TCR-pMHC complexes through repulsive particle interactions in previous models~\cite{DharanFarago2017}. We also do not account for attractive interaction forces that have been used to model TCR aggregation into microclusters. That is, we ignore size-based segregation forces, and membrane mediated potential of mean force interactions~\cite{DharanFarago2017,Siokisetal2018}. Note, however, as employed in~\cite{DharanFarago2017}, ignoring such forces allows insight into their relative role in driving the accumulation of TCR-pMHC complexes in the central region of synapses.

\begin{figure}[!tbp]
  \centering
  \includegraphics[width=.5\textwidth]{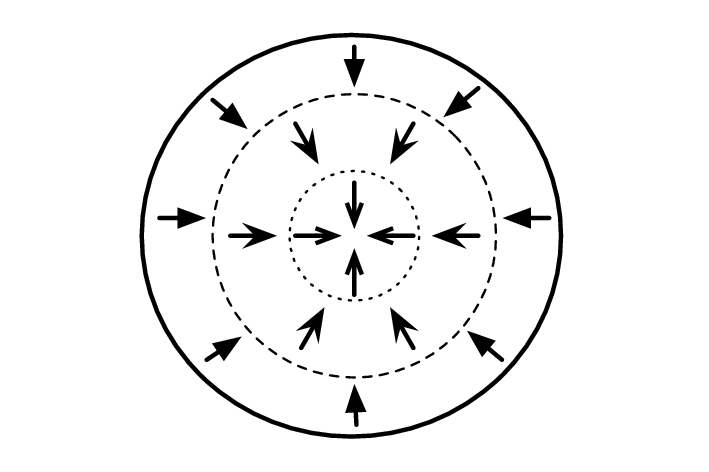}
  \caption{A graphical illustration of the contact area between the membranes of the T cell and the APC (antigen-presenting cell) as described in this section and also in \cite{DharanFarago2017,Siokisetal2018}. The dotted circle represents the actin-depletion region and has a radius of $\RC = 2\mu$m. In this region it is thought that dyenin motor-driven flows drive active, centripetal transport of TCR-pMHC complexes~\cite{DharanFarago2017}. Outside this region retrograde actin flows are thought to dominate active, centripetal transport of TCR-pMHC complexes~\cite{DharanFarago2017}. The dashed circle marks the edge of the pSMAC and has a radius of $\RP = 4\mu$m. Outside this region, centripetal actin retrograde flows are assumed to increase in strength~\cite{DharanFarago2017}. The arrows represent the active cytoskeleton forces given by \eqref{eq:ISpotential}.}
    \label{fig:ISForce}
\end{figure}

 The simulation domain is a circle centered at the origin with radius $R = 5.6\mu$m (solid circle in Fig.~\ref{fig:ISForce}), discretized into 13890 polygons, equipped with no-flux Neumann boundary conditions on the boundary of the domain. Model parameters and components are primarily taken from the model of~\cite{DharanFarago2017}, which only involved spatial transport, or the model of~\cite{Siokisetal2018}, which also included chemical reactions. As described in \cite{DharanFarago2017} the contact area between the T cell and the APC requires two types of adhesion bonds located in two different regions; therefore, the domain is divided into two regions. The inner circular region with a radius of $\RC = 2 \mu$m is called the actin-depletion region (dotted circle in Fig.~\ref{fig:ISForce})~\cite{DharanFarago2017}. The outer circular region with a radius of $\RP = 4\mu$m (dashed circle in Fig.~\ref{fig:ISForce}) represents the edge of the pSMAC (peripheral supra-molecular activation centers)~\cite{DharanFarago2017}. The system is initialized with 10 free TCRs and 10 pMHCs per polygon outside the actin-depletion region. The number of free TCRs and pMHCs are chosen to be consistent with the particle concentrations used in \cite{DharanFarago2017}. Following the models in~\cite{Siokis2017,Siokisetal2018}, we assume that free TCRs and pMHCs do not experience any active forces, and simply diffuse with $D = 0.1 \mu$m$^2$s$^{-1}$, while TCR-pMHC complexes diffuse with $D = 0.06\mu$m$^2$s$^{-1}$ and experience active forces that impart drift. The binding rate is set to $\lambda = 2\times 10^5 / 6.023$s$^{-1}$ and the binding radius is $\rb = 0.015 \mu$m, which corresponds to the average length of a TCR-pMHC complex. The TCR-pMHC complex unbinds at a rate of $\mu = 0.1$s$^{-1}$ \cite{Siokisetal2018}.

\begin{figure}[!tbp]
	\centering
	\includegraphics[width=1.05\textwidth]{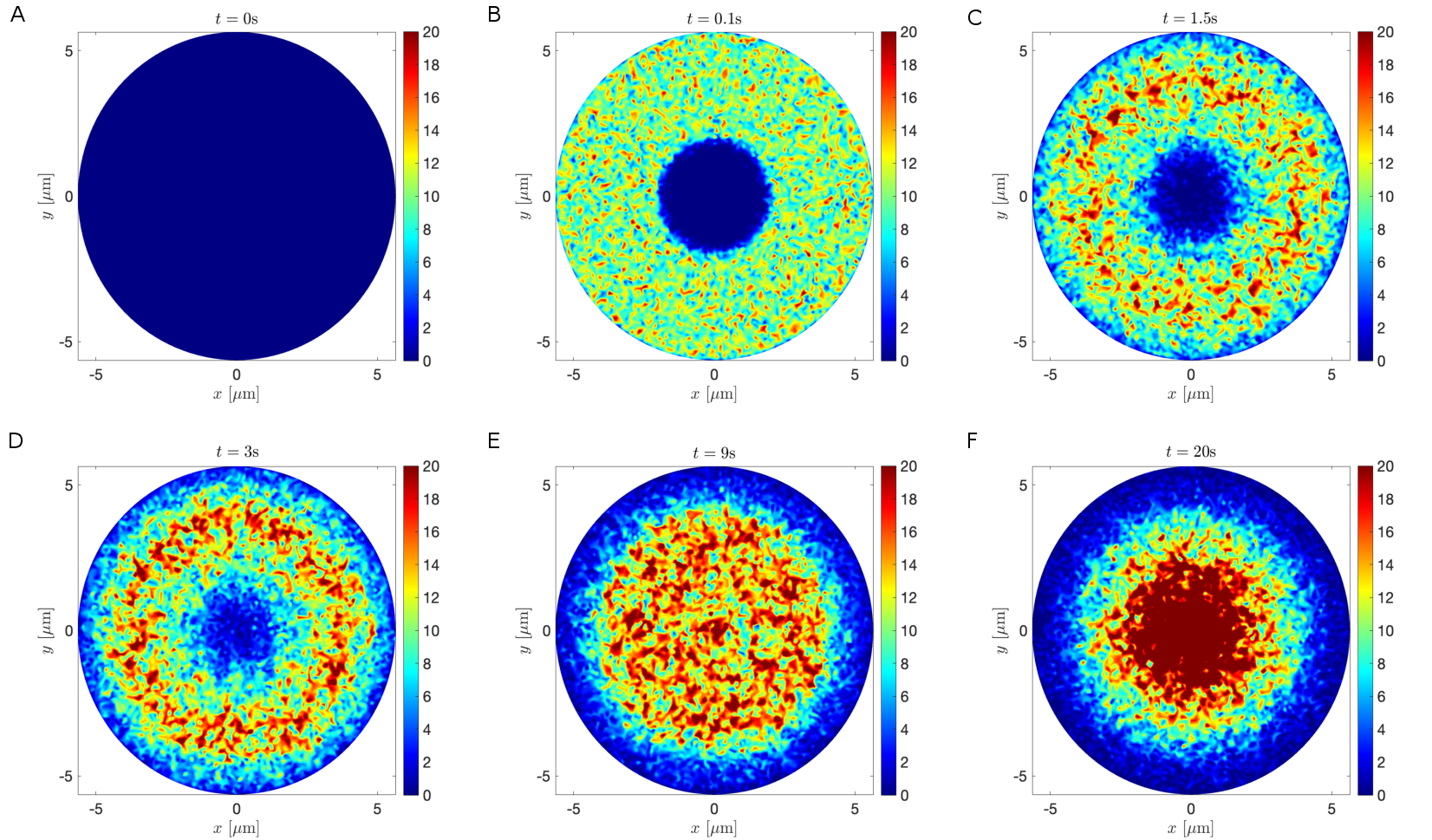}
	\caption{TCR-pMHC complex profiles at $t = 0$s, $0.1$s, $1.5$s, $3$s, $9$s, and $20$s from one CRDDME simulation. We observe that TCR-pMHC complexes accumulate within the central region over timescales of seconds in the absence of attractive two-body interactions as was also observed in~\cite{DharanFarago2017}.}
	  \label{fig:ISTimeEvo}
  \end{figure}
Once a free TCR binds with a free pMHC, the TCR-pMHC complex is then actively transported by dynein motors along microtubules, and via an actin retrograde flow, which together induce a centripetal force that drives the TCR-pMHC complexes to the central region of the domain~\cite{DharanFarago2017}. See Fig.~\ref{fig:ISForce} for details on where each force dominates. Modeling these two active processes, especially the dynein-induced forces, is highly challenging as the motor activity depends on various factors such as free ATP and the motor density. To bypass such complexity, we followed the approach used in \cite{DharanFarago2017,Siokis2017,Siokisetal2018} and modeled the impact of the retrograde actin flow via a centripetal-induced force via the potential function used in~\cite{DharanFarago2017}
\begin{equation}
	\phi(x,y) = \begin{cases}
		f_0\sqrt{x^2 + y^2}, &\sqrt{x^2 + y^2} \leq \RP, \\
		2f_0\paren{\sqrt{x^2 + y^2} - \RP} + f_0\RP, &\sqrt{x^2 + y^2} > \RP.
	\end{cases}
	\label{eq:ISpotential}
\end{equation}
Here $f_0 = 1$ corresponds to the magnitude of the TCR-pMHC centrally directed force in~\cite{Siokisetal2018}. Note, the strength of the force increases as one moves outward from the pSMAC boundary, modeling increased centripetal actin retrograde flows that have been observed in that region~\cite{DharanFarago2017}. 

In Fig.~\ref{fig:ISTimeEvo} we demonstrate TCR-pMHC complex profiles at $t = 0$s, $0.1$s, $1.5$s, $3$s, $9$s, and $20$s from one CRDDME SSA simulation. We observe that our CRDDME-based model can reproduce the in silico knock-down assays of~\cite{DharanFarago2017} (Fig. 4) in which attractive two-body interactions were removed. Comparing to Fig. 4 in \cite{DharanFarago2017}, we observe similar results, with TCR-pMHC complexes accumulating within the central region over a timescale of seconds. Since attractive two-body membrane-induced forces were not included, we do not observe microclusters of TCR-pMHC complexes as in the full models of~\cite{Siokisetal2018,DharanFarago2017}. These microclusters will move more slowly towards the central core due to their larger size and repulsive interactions with other synapse components, which was found to drive the switch from TCR-pMHC clusters accumulating in the central region within seconds to more physiological timescales of tens of minutes (as observed when such forces are included within the models~\cite{DharanFarago2017,Siokisetal2018}).

\section{Discussion}
\label{sec:discussion}
By using a hybrid discretization method combining the edge-averaged finite element discretization of the drift-diffusion term and a finite volume discretization of reaction terms, we developed a convergent unstructured mesh jump process approximation, the convergent-reaction-drift-diffusion master equation (CRDDME), to the abstract volume-reactivity model with drift for reversible reactions. The CRDDME is capable of handling general bimolecular reaction kernels when the motion of molecules is influenced by background one-body potential fields, supports unstructured polygonal meshes in two or three dimensions, and should be extendable to surfaces~\cite{MaIsaacsonSurfCRDME2024}. One benefit to the CRDDME approach is that it is consistent with key physical properties of the corresponding volume reactivity model. In particular, the reaction rates are consistent with detailed balance of reaction fluxes holding at equilibrium, and in the absence of reactions discrete versions of equilibrium Gibbs-Boltzmann distributions are recovered. Moreover, when modeling multiparticle systems, reactive transition rates are naturally given in terms of the rates derived for a minimal two-particle, $\textrm{A} + \textrm{B} \to \textrm{C}$, and minimal one-particle, $\textrm{C} \to \textrm{A} + \textrm{B}$, reaction.

Another advantage to the CRDDME approach is that it will be consistent with the RDME model in its treatment of diffusive transport and linear reactions, yet converges to the underlying spatially continuous volume-reactivity model. This enables the CRDDME model to potentially reuse many RDME-based extensions for modeling spatial transport. We note that the spatial drift-diffusion transition rates we obtained could also be used without modification to define unstructured-grid RDME models that support drift due to one-body potentials.

There are still a number of directions in which the CRDDME could be improved. An interesting future direction would be to consider time-dependent potentials. Such a generalization would enable the study of time-dependent actin flows in T cell signaling. Extending the CRDDME to support two-body interaction potentials, as begun in~\cite{HeldmanTheis2023}, would be useful to model volume exclusion in dense particle systems and other types of particle interactions~\cite{Siokisetal2018,mintonCrowding2001,EngblomEtAl2018}. As discussed in the previous section, such features have been found to be important to accurately capture aggregation dynamics during T cell synapse formation~\cite{Siokisetal2018,DharanFarago2017}. This work illustrates empirical convergence of the CRDDME to the VR model for several statistics. It would also be of interest to prove convergence of the hybrid discretization-based CRDME/CRDDME to the VR model, as begun in~\cite{HeldmanTheis2023,Heldman2024}, along with examining stronger types of convergence of the jump process associated with the CRDDME to the process associated with the VR model. Finally, we note that the use of mass lumping for advection-dominated problems has been shown to induce dispersion errors in finite element methods, and a number of studies have investigated corrections to help minimize such errors~\cite{WendlandSchulz2005,GuermondPasquetti2013}. Whether such methods can be adapted to derive CRDDME-type models, i.e. give discretizations that correspond to transition rate matrices, is currently an open problem.

\begin{acknowledgements}
	SAI and YZ were supported by National Science Foundation DMS-1902854.
\end{acknowledgements}


\bibliographystyle{apsrev4-1.bst}

\bibliography{lib.bib}
\clearpage
\begin{appendix}

\section{Accuracy of the EAFE Method}
\label{app:EAFEAccuracy}
In this section, we demonstrate the second order accuracy of the EAFE discretization method. For all examples in this section, we numerically solve a time-independent PDE using either the spatial operator for the density from \eqref{eq:semidiscreteFPE}, or the effective transition rate matrix from \eqref{eq:discreteDDOperator}. We choose $\Omega$ to be either a square domain or a circle.

As a test problem, we consider the following steady-state problem on a square $[-0.5, 1.5] \times [-0.5, 1.5]$ and on the circle centered at $(-0.5, 0)$ with radius one with a no-flux boundary condition
 \begin{equation}
	\begin{aligned}
		-\nabla\cdot D\brac{\nabla\rho(\vx) + \rho(\vx)\nabla\phi(\vx)} + \rho(\vx) &= f(\vx), &\vx \in \Omega \\
  		D\brac{\nabla\rho(\vx) + \rho(\vx)\nabla\phi(\vx)}\cdot\eta(\vx) &= 0, &\vx \in \partial\Omega
  	\end{aligned}
  	\label{eq:TestExample}
\end{equation}
where the function describing the potential field is given by
\begin{equation}
	\phi(\vx) = x_1^2 + x_2^2.
	\label{eq:potentFcn}
\end{equation} For both cases, we assume a constant diffusion constant, with $D = 1$ inside the square and $D = 10$ inside the circle. On the square, the forcing term was taken to be
\begin{align*}
	f(\vx) = -e^{-x_1^2 - x_2^2}[(-1-8\pi^2)&\cos(2\pi x_1)\cos(2\pi x_2) \\
			&+4\pi(x_1\cos(2\pi x_2)\sin(2\pi x_1) + x_2\cos(2\pi x_1)\sin(2\pi x_2))],
\end{align*}
while on the circle we chose
\begin{equation*}
	f(\vx) = e^{-x_1^2 - x_2^2}.
\end{equation*}
These choices correspond to analytical solutions given by
\begin{equation*}
	\rho(\vx) = e^{-x_1^2-x_2^2}\cos(2\pi x_1)\cos(2\pi x_2)
\end{equation*}
on the square and
\begin{equation*}
	\rho(\vx) = e^{-x_1^2-x_2^2}
\end{equation*}
on the circle.

We triangulate the domain and approximate $\rho(\vx)$ by
\begin{equation*}
	\rho_h(\vx) = \sum_{i = 1}^N \rho_i \psi_{i}(\vx),
\end{equation*}
with $\rho_i \approx \rho(\vx_i)$, $i = 1, \dots, N$, the nodal solution values on our triangulation. Let $\vec{\rho}_h = \brac{\rho_1, \dots, \rho_N}^T$. The EAFE discretization of \eqref{eq:TestExample} gives
\begin{equation}
	(-S + M)\vec{\rho}_h = \vec{F},
\end{equation}
with $S$ the EAFE stiffness matrix, and $M$ the corresponding lumped mass matrix we previously derived. Here $\vec{F}$ is a vector of forcing terms given by
\begin{equation*}
	\vec{F}_i = \int_{\Omega}\psi_i(\vx)f(\vx)\, d\vx.
\end{equation*}
In solving the steady-state problem \eqref{eq:TestExample}, this forcing vector is evaluated using numerical quadrature.

\begin{figure}[!tbp]
  \centering
  \includegraphics[width = 1\textwidth]{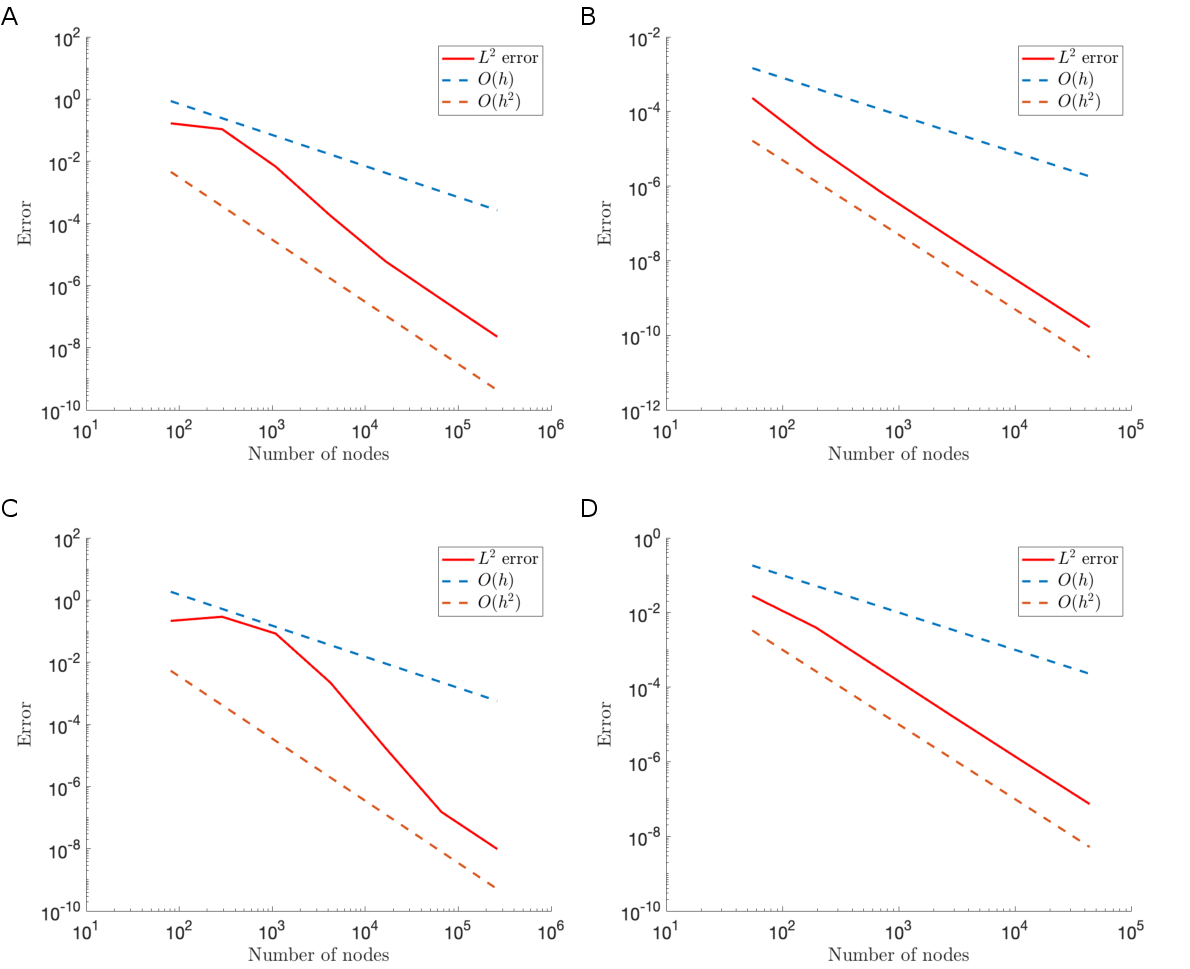}
  \caption{Convergence of the EAFE approximation to the true solution using the potential given by \eqref{eq:potentFcn} and \eqref{eq:potentFcnLarge} as the triangular mesh is refined. In panel A and C, we plot error of the EAFE approximation on the square using \eqref{eq:potentFcn} and \eqref{eq:potentFcnLarge} respectively. In panel B and D, we plot the error of the EAFE approximation on the circle using \eqref{eq:potentFcn} and \eqref{eq:potentFcnLarge} respectively.}
  \label{fig:EAFE_ConvergTests}
\end{figure}

In Fig.~\ref{fig:EAFE_ConvergTests}A and B, we plot the error between the discerte solution and the analytical solution as the mesh size is successively halved. In both cases, we observe a second-order convergence.
To examine the convergence in a convection-dominated scenario, we increase the magnitude of the potential \eqref{eq:potentFcn} by a factor of $30$ to
\begin{equation}
	\phi(\vx) = 30(x_1^2 + x_2^2).
	\label{eq:potentFcnLarge}
\end{equation}
In Figs.~\ref{fig:EAFE_ConvergTests}C and D, we plot the error between the discrete solution and the PDE solution as the mesh size is successively halved. We again observe a consistent second-order convergence rate on the circle, with the convergence rate on the square becoming second order once the mesh is sufficiently small.
\begin{figure}[!tbp]
  \centering
  \includegraphics[width = 1\textwidth]{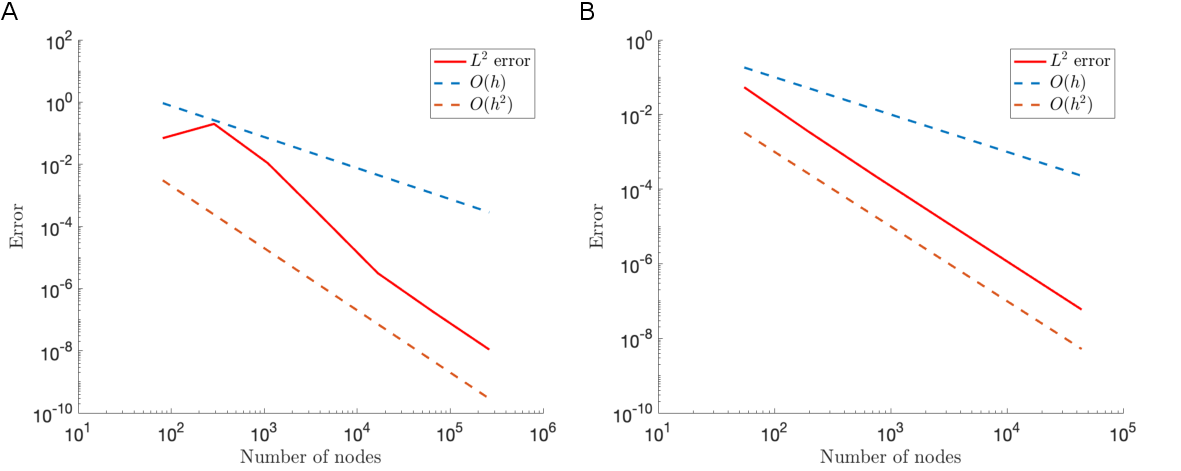}
  \caption{Convergence of the EAFE approximation to the true solution using the potential given by \eqref{eq:potentFcnTwoWell} as the mesh is refined. In panel A, we plot error of the EAFE approximation on a square of side length $1$. In panel B, we plot the error of the EAFE approximation on a circle centered at the origin of radius $1$.}
  \label{fig:EAFE_Converg_TwoWell}
\end{figure}
Finally, we repeat the same analysis on both the square and the circle using a two-well potential similar to the one in \cite{Schutteetal2011,VandenEijnden2010}
\begin{equation}
	\phi(\vx) = \frac{5}{2}(1-x_1^2)^2 + 5x_2^2.
	\label{eq:potentFcnTwoWell}
\end{equation}
For both the square and circle domains, we again observe second-order convergence as the mesh size is decreased (Fig.~\ref{fig:EAFE_Converg_TwoWell}).

An empirical estimate for the rate of convergence of $\rho_h$ to $\rho$ in the $L^2$ norm can be calculated using the approach of \cite{GriffithPeskin2005}. Denote by $e_2^N$ the $L^2$ norm of the difference in the two approximations obtained using a mesh of $N^2$ nodes and $(2N)^2$ nodes,
\begin{equation}
	e_2^N = \norm{\rho_h^N - \mathcal{I}^{2N\to N}\rho_h^{2N}}_2,
\end{equation}
where $\mathcal{I}^{2N\to N}$ is the interpolation from the finer to the coarser meshes. We note that in our approach of mesh refinement via dividing each triangle into four congruent triangles \cite{Carey1997}, $\mathcal{I}^{2N\to N}$ is simply the $N^2\times N^2$ identity matrix as the finer mesh is obtained by subdividing each triangle into four congruent triangles; and therefore the the first $N^2$ mesh points in the finer mesh correspond to the $N^2$ mesh points in the coarser mesh. The empirical estimate for the rate of convergence is then defined to be
\begin{equation}
	r_2^N = \log_2\paren{\frac{e_2^N}{e_2^{2N}}}.
	\label{eq:convergRateMesh}
\end{equation}
Table.~\ref{tab:empiricalConvergRate} summarizes the empirical convergence rates in $L^2$ norm computed using \eqref{eq:convergRateMesh} for the approximated solution to \eqref{eq:TestExample} on both the square domain and the circle domain. In all three cases, we observe a rate of convergence that is at least second order when solving \eqref{eq:TestExample} on a square and a circle.
\begin{table}[!h]
\begin{tabular}{|c|c|c|c|c|c|c|}
\hline
\multicolumn{1}{|c|}{ } & \multicolumn{2}{c|}{$\phi(\vx) = x_1^2 + x_2^2$} & \multicolumn{2}{c|}{$\phi(\vx) = 30x_1^2 + 30x_2^2$} & \multicolumn{2}{c|}{$\phi(\vx) = \frac{5}{2}(1-x_1^2)^2 + 5x_2^2$} \\ \hline
$N$ & $r_2$, square & $r_2$, circle & $r_2$, square & $r_2$, circle & $r_2$, square & $r_2$, circle \\ \hline
194 &0.92386 &1.5813  &2.2812  &3.1383  &0.4107 &2.2258\\ \hline
727 &1.8901 &1.9203  &0.19782  &1.7173  &2.0042 &2.0509\\ \hline
2813 &2.0127 &2.0224  &2.60487  &1.8654  &2.2029 &1.9724\\ \hline
11065 &2.1593 &2.0309  &2.90366  &2.1669 &2.1156 &1.9891\\ \hline
43889 &2.0848 & &2.763210  & &2.0126 &\\ \hline
\end{tabular}
\caption{Empirical convergence rates for the approximated solution to \eqref{eq:TestExample} using three different potential functions in the $L^2$ norm}.
\label{tab:empiricalConvergRate}
\end{table}

\section{Multi-particle $\textrm{A} + \textrm{B} \leftrightarrows \textrm{C}$ reaction}
\label{app:multipartModel}
In this section we consider the multiparticle $\textrm{A} + \textrm{B} \leftrightarrows \textrm{C}$ reaction within a bounded domain $\Omega \subset \R^d$, $d \leq 3$. We will show that the discretization used in Sections \ref{sec:DriftDiffApprox} and \ref{sec:FVM_RXs} to derive the master equation approximation for two particles \eqref{eq:twopartCRDDME} naturally leads to a master equation approximation for the multiparticle system, with transition rates for multiparticle reactions given in terms of transition rates for the two-particle case. The final rates we derive are summarized in Table~\ref{Tab:CRDDMETransitionRates}.

We begin by formulating the multiparticle abstract volume reactivity model with $\textrm{A} + \textrm{B} \leftrightarrows \textrm{C}$ reaction, using a similar notation to \cite{IsaacsonCRDME2013,IsaacsonZhang17}. Let $A(t)$ denote the stochastic process for the number of \textrm{A} molecules in the system at time $t$, with $B(t)$ and $C(t)$ defined similarly. Values of $A(t)$, $B(t)$, and $C(t)$ are given by $a$, $b$, and $c$ respectively (i.e. $A(t) = a$). Denote by $\vQa_l(t) \in \Omega$ the stochastic process for the position of the $l^\text{th}$ \textrm{A} molecule within the domain $\Omega$, with $\vqa_l$ be a possible value of $\vQa(t)$. The position vector of \textrm{A} molecules when $A(t) = a$ is given by
\begin{equation*}
  \vQa(t) = (\vQa_1(t), \dots, \vQa_a(t)) \in \Omega^{a},
\end{equation*}
where $\Omega^a = \Omega \times \cdots \times \Omega \subset \R^{ad}$. In a similar way, $\vqa$ denotes a possible value of $\vQa(t)$,
\begin{equation*}
  \vqa = \vQa(t) = (\vqa_1, \dots, \vqa_a).
\end{equation*}
$\vQb_m(t)$, $\vQb(t)$, $\vQc_n(t)$, $\vQc(t)$, $\vqb_m$, $\vqb$, $\vqc_n$, $\vqc$, $\Omega^b$, and $\Omega^c$ are defined analogously. We denote by $\Phi^a(\vqa)$ the background one-body potential for molecules of species \textrm{A} when $A(t) = a$ with $\vQa(t) = \vqa$. In what follows, we assume particles of species \textrm{A} experience the one-body potential, $\phi^{\textrm{A}}(\vx)$, so that when the species \textrm{A} state vector is $\vqa$ it experiences the collective potential
\begin{equation*}
	\begin{aligned}
		\Phi^a(\vqa) &:= \sum_{\ell=1}^a \phi^{\textrm{A}}(\vqa_{\ell}).
	\end{aligned}
\end{equation*}
$\Phi^b(\vqb)$, $\Phi^c(\vqc)$, $\vec{\sigma}^b$, and $\vec{\sigma}^c$ are defined similarly.

With this notation, we let $\fabc(\vqa,\vqb,\vqc,t)$ denote the probability density of having $A(t) = a$, $B(t) = b$, and $C(t) = c$ with $\vQa(t) = \vqa$, $\vQb(t) = \vqb$, and $\vQc(t) = \vqc$. Here we assume that molecules of the same species are indistinguishable and that $\fabc(\vqa, \vqb, \vqc)$ is symmetric respect to permutations of the molecule orderings within $\vqa$, $\vqb$, and $\vqc$. With this assumption, $\fabc$ is normalized so that
\begin{equation*}
  \sum_{a=0}^{\infty} \sum_{b=0}^{\infty} \sum_{c=0}^{\infty} \brac{ \frac{1}{a!\, b!\, c!}
  \int_{\Omega^{a}} \int_{\Omega^{b}} \int_{\Omega^{c}} \fabc\paren{\vqa,\vqb,\vqc,t} \, d\vqc \, d\vqb \, d\vqa
  } = 1.
\end{equation*}
Finally, we denote by $\vf(t)$ the overall probability density vector at time $t$
\begin{equation*}
  \vf(t) = \{ f^{(a,b,c)}(\vqa,\vqb,\vqc,t)\}_{a,b,c}.
\end{equation*}
This density function then satisfies the forward Kolmogorov equation
\begin{equation}
	\label{eq:multipartABtoCEqs}
	\PD{\vf}{t}(t) = (\diffop + \Rp + \Rm)\vf(t).
\end{equation}
Here we assume molecules cannot leave the domain $\Omega$, and as such a reflecting Neumann boundary condition is imposed for each component $f^{(a,b,c)}(\vqa,\vqb,\vqc,t)$ on the domain boundary $\partial \paren{\Omega^{a+b+c}}$. The linear operators $\diffop$, $\Rp$, $\Rm$ correspond to drift-diffusion, the association reaction, and the dissociation reaction respectively. The drift-diffusion operator is given by
\begin{equation}
  \label{eq:multipartDiffop}
  \begin{aligned}
  	(\diffop \vf(t))_{a,b,c} =
  	&\sum_{l = 1}^a\nabla_{\vqa_l}\cdot\brac{\DA\paren{\nabla_{\vqa_l}f^{(a,b,c)} + f^{(a,b,c)}\nabla_{\vqa_l}\phi^{\textrm{A}}(\vqa)}} \\
  	&+\sum_{m = 1}^b\nabla_{\vqb_m}\cdot\brac{\DB\paren{\nabla_{\vqb_m}f^{(a,b,c)} + f^{(a,b,c)}\nabla_{\vqb_m}\phi^{\textrm{B}}(\vqb)}} \\
  	&+\sum_{n = 1}^c\nabla_{\vqc_n}\cdot\brac{\DC\paren{\nabla_{\vqc_n}f^{(a,b,c)} + f^{(a,b,c)}\nabla_{\vqc_n}\phi^{\textrm{C}}(\vqc)}},
  \end{aligned}
\end{equation}
where $\nabla_{\vqa_l}$ denotes the $d-$dimensional gradient operator with respect to $\vqa_l$, and $\nabla_{\vqb_m}$ and $\nabla_{\vqc_n}$ are defined similarly. Notice that this operator naturally splits into transport operators for each individual particle's motion, so that in the absence of reactions particles would move by independent drift-diffusion processes. Moreover, this splitting means we can immediately apply our single-particle EAFE approximation to each operator independently to derive transition rates for each particle.

To define the reaction operators, $\Rp$ and $\Rm$, we make use of the following notations for adding an \textrm{A} molecule at $\vx$ or removing the $l^{\textrm{th}}$ molecule of species \textrm{A} from a given state, $\vqa$,
\begin{align*}
	\vqa \cup \vx &= (\vqa_1, \dots, \vqa_a, \vx), & \vqa \setminus \vqa_l &= (\vqa_1, \dots, \vqa_{l-1}, \vqa_{l+1}, \dots, \vqa_a).
\end{align*}
With these definitions, the reaction operator for the $\textrm{A} + \textrm{B} \to \textrm{C}$ association reaction is given by
\begin{equation}
\label{eq:multipartRp}
  \begin{aligned}
    (\Rp \vf(t))_{a,b,c} =
    &-\paren{\sum_{l=1}^a \sum_{m=1}^b \kp \paren{\vqa_{l}, \vqb_{m}}} \fabc(\vqa,\vqb,\vqc,t) \\
    &+ \sum_{n=1}^{c} \!\brac{\int_{\Omega^2} \kp(\vqc_n \vert \vx, \vy)
      f^{(a+1,b+1,c-1)}(\vqa \cup \vx, \vqb \cup \vy, \vqc \setminus \vqc_n, t) d\vx d\vy}\!\!.
  \end{aligned}
\end{equation}
The reaction operator for the $\textrm{C} \to \textrm{A} + \textrm{B}$ dissociation reaction is given by
\begin{multline} \label{eq:multipartRm}
  (\Rm \vf(t))_{a,b,c} = - \paren{ \sum_{n=1}^{c} \km(\vqc_n)} \fabc(\vqa, \vqb, \vqc, t) \\
  + \sum_{l=1}^a \sum_{m=1}^b \brac{\int_{\Omega} \km\paren{\vqa_l,\vqb_m\vert\vz}
    f^{(a-1,b-1,c+1)}\paren{\vqa\setminus \vqa_l, \vqb\setminus \vqb_m, \vqc \cup \vz, t} d\vz}.
\end{multline}

To arrive at a master equation model approximation of \eqref{eq:multipartABtoCEqs}, we utilize the discretizations for drift-diffusion in Section \ref{sec:DriftDiffApprox} and reactions developed in Section \ref{sec:FVM_RXs}. We again consider a polygonal mesh approximation of the domain $\Omega$, denoted by $\{V_i\}_{i = 1}^N$. We let $\bi^a = (i^a_1,\dots,i^a_a)$ denote the multi-index label of the hyper-voxel $\vV_{\bi^a}$
\begin{equation*}
  \vV_{\bi^a} := V_{i^a_1} \times \dots \times V_{i^a_a},
\end{equation*}
with $\bj^b$, $\bk^c$, $\vV_{\bj^b}$, and $\vV_{\bk^c}$ defined similarly. We will use $\vV_{\bi^a \bj^b} := \vV_{\bi^a} \times \vV_{\bj^b}$ and
$\vVabc := \vV_{\bi^a} \times \vV_{\bj^b} \times \vV_{\bk^c}$ to denote the multi-species hypervoxels. A well-mixed piecewise constant approximation of the probability the system is in the state $A(t) = a$, $B(t) = b$, $C(t) = c$ with $(\vQa(t), \vQb(t), \vQc(t)) \in \vVabc$ at time $t$ is given by
\begin{subequations}
  \begin{align}
    \Fabcijk(t) &:= \prob \brac{A(t) = a, B(t) = b, C(t) = c \text{ and }
        \paren{\vQa(t), \vQb(t), \vQc(t)} \in \vVabc} \notag \\
        &\phantom{:}= \frac{1}{a! b! c!} \int_{\vVabc} \fabc(\vqa,\vqb,\vqc,t) \, d\vqc \, d\vqb \, d\vqa \notag \\
        &\phantom{:}\approx \frac{\abs{\vVabc}}{a! \, b! \, c!} \fabc(\vqa_{\bi^a}, \vqb_{\bj^b}, \vqc_{\bk^c},t),
        \label{eq:FabcijkDef}
  \end{align}
\end{subequations}
where $\vqa_{\bi^a}$ denotes the vector containing the centroids of the voxels in $V_{\bi^a}$, with $\vqb_{\bj^b}$ and $\vqc_{\bk^c}$ defined similarly. We denote by $\vF(t)$ the overall state vector
\begin{equation*}
  \vF(t) = \{ \Fabcijk(t) \}_{\bi^a,\bj^b,\bk^c}.
\end{equation*}
We follow a similar approach as we employed in~\cite{IsaacsonZhang17}. Using the EAFE discretization as in Section \ref{sec:DriftDiffApprox} to approximate each single-particle drift-diffusion operator within $\diffop$, integrating the action of each of the reaction operators, $\Rp$ and $\Rm$, on $\fabc$ over $\vVabc$, and applying the locally well-mixed (i.e. piecewise constant) approximation \eqref{eq:FabcijkDef}, the state vector $\vF(t)$ satisfies the master equation
\begin{equation}
	\D{\vF}{t}(t) = \paren{\diffoph + \Rph + \Rmh} \vF(t).
	\label{eq:masterEqnF}
\end{equation}
The discretized drift-diffusion operator $\diffoph$ is given by
\begin{equation}
	\begin{aligned}
  (\diffoph \vF(t))_{\bi^a,\bj^b,\bk^c}
  &= \sum_{l=1}^{a} \sum_{i'= 1}^{N} \brac{ \paren{\diffophA}_{i^a_l i'} F_{\bi^a\setminus i^a_l\cup i', \bj^b, \bk^c}(t)
    - \paren{\diffophA}_{i' i^a_l } \Fabcijk(t)} \\
  &+ \sum_{m=1}^{b} \sum_{j'= 1}^{N} \brac{ \paren{\diffophB}_{j^b_m j'} F_{\bi^a, \bj^b\setminus j^b_m\cup j', \bk^c}(t)
    - \paren{\diffophB}_{j' j^b_m } \!\Fabcijk(t)} \\
  & + \sum_{n=1}^{c} \sum_{k'= 1}^{N} \brac{ \paren{\diffophC}_{k^c_n k'} F_{\bi^a,\bj^b,\bk^c\setminus k^c_n\cup k'}(t)
    - \paren{\diffophC}_{k' k^c_n } \!\Fabcijk(t)},
	\end{aligned}
\end{equation}
where $\diffophA$, $\diffophB$, and $\diffophC$ are the single-particle species-specific EAFE discretized drift-diffusion operators, with potentials given by $\phi^{\textrm{A}}$, $\phi^{\textrm{B}}$, and $\phi^{\textrm{C}}$ respectively.

The association reaction operator $\Rph$ is given by
\begin{multline}
  (\Rph \vF(t))_{\bi^a,\bj^b,\bk^c} =
  -\paren{\sum_{l=1}^a \sum_{m=1}^b \kp_{i^a_l j^b_m} } \Fabcijk(t) \\
  + \frac{(a+1)(b+1)}{c} \sum_{n=1}^{c} \sum_{i'=1}^{N} \sum_{j'=1}^{N}
  \kp_{i' j' k_n^c} F_{\bi^a \cup i', \bj^b \cup j', \bk^c \setminus k_n^c}(t),
\end{multline}
and the backward reaction dissociation operator $\Rmh$ is given by
\begin{equation}
  (\Rmh \vF(t))_{\bi^a,\bj^b,\bk^c} =
  - \paren{ \sum_{n=1}^{c} \km_{k_n^c}} \Fabcijk(t)
  + \frac{c+1}{a b} \sum_{l=1}^a \sum_{m=1}^b \sum_{k'=1}^N \km_{i_l^a j_m^b k'}
  F_{\bi^a \setminus i_l^a, \bj^b \setminus j_m^b, \bk^c \cup k'}(t).
\end{equation}

Denote by $A_i(t)$ the stochastic process for the number of molecules of species $\textrm{A}$ in voxel $V_i$ at time $t$, with $B_j(t)$ and $C_k(t)$ defined similarly. For each species, the vector of the state at time $t$ is given by
 \begin{equation*}
 	\vA(t) = (A_1(t), A_2(t), \dots, A_N(t)),
 \end{equation*}
 with $\vB(t)$ and $\vC(t)$ defined similarly. The vector $\va$ denotes the corresponding values of $\vA(t)$ and is given by
 \begin{equation*}
 	\va = (a_1, a_2, \dots, a_N).
 \end{equation*}
 The vectors $\vb$ and $\vc$ are defined similarly. Finally, we let
 \begin{equation*}
 	P(\va, \vb, \vc, t) := \prob\brac{\vA(t) = \va, \vB(t) = \vb, \vC(t) = \vc}.
 \end{equation*}

 Following a similar analysis as in \cite{IsaacsonZhang17}, we can derive a master equation for $\{P(\va, \vb, \vc, t)\}_{\va, \vb, \vc}$
 \begin{equation}
 	\D{\vP}{t}(t) = \paren{\diffoph + \Rph + \Rmh} \vP(t)
 	\label{S1eq:masterEqnProb}
 \end{equation}
 from the master equation for $\vF(t)$. Here the drift-diffusion operator is then given by
 \begin{equation*}
 	\begin{aligned}
  	(\diffoph \vP)(\va, \vb, \vc, t)
  	&= \sum_{i=1}^{N} \sum_{i'= 1}^{N} \brac{ \paren{\diffophA}_{i i'} (a_{i'}+1)P(\va+\ve_{i'}-\ve_i, \vb, \vc, t)
    	- \paren{\diffophA}_{i' i} P(\va, \vb, \vc, t)} \\
  	&+ \sum_{j=1}^{N} \sum_{j'= 1}^{N} \brac{ \paren{\diffophB}_{j j'} (b_{j'}+1)P(\va, \vb+\ve_{j'}-\ve_j, \vc, t)
    	- \paren{\diffophB}_{j' j} P(\va, \vb, \vc, t)} \\
  	& + \sum_{k=1}^{N} \sum_{k'= 1}^{N} \brac{ \paren{\diffophC}_{k k'} (c_{k'}+1) P(\va, \vb, \vc + \ve_{k'} - \ve_k, \vc, t)
    	- \paren{\diffophC}_{k' k} P(\va, \vb, \vc, t)},
 	\end{aligned}
 	\label{S1eq:diffophProb}
 \end{equation*}
 where $\ve_i$ denotes the unit vector along the $i$th coordinate axis of $\R^N$.

 Similarly, the association operator is given by
 \begin{equation*}
 	\begin{aligned}
  	(\Rph \vP)(\va, \vb, \vc, t) =
  	&-\paren{\sum_{i=1}^N \sum_{j=1}^N \kp_{i j}a_ib_j} P(\va, \vb, \vc, t) \\
  	&+ \sum_{i=1}^{N} \sum_{j=1}^{N} \sum_{k=1}^{N}
  	\kp_{i j k}(a_i+1)(b_j+1)P(\va+\ve_i, \vb+\ve_j, \vc-\ve_k, t),
	\end{aligned}
	\label{S1eq:RphProb}
 \end{equation*}
and the backward reaction dissociation operator $\Rmh$ is given by
\begin{equation*}
	\begin{aligned}
  	(\Rmh \vP)(\va, \vb, \vc, t) =
  	&- \paren{ \sum_{k=1}^{N} \kmk c_k} P(\va, \vb, \vc, t) \\
  	&+ \sum_{i=1}^N \sum_{j=1}^N \sum_{k=1}^N \kmijk (c_k+1)P(\va-\ve_i, \vb-\ve_j, \vc+\ve_k, t).
	\end{aligned}
	\label{S1eq:RmhProb}
\end{equation*}
We find that the effective transition rates for the mulitparticle system are then as given in Table~\ref{Tab:CRDDMETransitionRates}, fully determined by the single-particle transport rates, the two particle association rate $\kp_{i j k}$, and the one-particle dissociation rate $\km_{i j k}$.

\section{Pointwise detailed balance for multi-particle $\textrm{A} + \textrm{B} \leftrightarrows \textrm{C}$ reaction}
\label{app:multipartDB}
 Denote by $\bar{\vP}$ the equilibrium solution of \eqref{S1eq:masterEqnProb}. We expect it to satisfy the detailed balance condition that
 \begin{equation}
 	a_ib_j\kpijk \bar{\vP}(\va, \vb, \vc) = (c_k+1)\kmijk \bar{\vP}(\va-\ve_i, \vb-\ve_j, \vc+\ve_k).
 	\label{eq:DBMultPartPabc}
 \end{equation}
 As derived in \cite{IsaacsonZhang17}, we have that
 \begin{equation}
 	\bar{\vP}(\va, \vb, \vc) = \frac{a!b!c!}{\prod_{i = 1}^N a_i! b_i! c_i!}\bar{\vF}_{\bi^a\bj^b\bk^c},
 	\label{eq:FabcijkPabc}
 \end{equation}
 where $\bar{\vF}_{\bi^a\bj^b\bk^c}$ denotes the steady state of \eqref{eq:masterEqnF}. Here we assume that the two representations of state are chosen to be consistent so that
 \begin{equation*}
 	a_i = \abs{\{i_l^a \vert i_l^a = i, l = 1, \dots, a\}},
 \end{equation*}
 with $\abs{\cdot}$ represents the cardinality of the set. $b_i$ and $c_i$ are defined similarly. Using \eqref{eq:DBMultPartPabc} and \eqref{eq:FabcijkPabc}, we then have the detailed balance condition for $\bar{\vF}_{\bi^a\bj^b\bk^c}$ given by
 \begin{equation}
 		ab\kpijk\bar{\vF}_{\bi^a\bj^b\bk^c} = (c+1)\kmijk\bar{\vF}_{\bi^a\setminus i, \bj^b\setminus j, \bk^c\cup k},
		\label{eq:correctDBMultPartF}
 \end{equation}
 which immediately implies that
 \begin{equation}
 	(\Rph \bar{\vF})_{\bi^a\bj^b\bk^c} = (\Rmh \bar{\vF})_{\bi^a\bj^b\bk^c}.
 	\label{eq:DBMultPartF}
 \end{equation}
 This implies that
 \begin{equation*}
 	(\diffoph \bar{\vF})_{\bi^a\bj^b\bk^c} = 0,
 \end{equation*}
 so that
 \begin{equation}
 	\bar{\vF}_{\bi^a \bj^b \bk^c} = \frac{\abs{\vVabc}}{(\hat{Z}_A)^a(\hat{Z}_B)^b(\hat{Z}_C)^c}e^{-\Phi^a_{\bi^a}-\Phi^b_{\bj^b}-\Phi^c_{\bk^c}}\hat{\pi}(a,b,c).
 	\label{eq:Fss}
 \end{equation}
 Here $\Phi^a_{\bi^a} = \sum_{\ell = 1}^a \phi_{i_{\ell}^a}^A$ denotes the total potential for species $A$ when the particles of species \textrm{A} are at the positions in $\bi^a$. The total potentials for species $B$ and $C$, given by $\Phi^b_{\bj^b}$ and $\Phi^c_{\bk^c}$, are defined similarly. The multi-dimensional normalization factor $(\hat{Z}_A)^a$ is given by
 \begin{equation*}
 	(\hat{Z}_A)^a = \brac{\sum_{j = 1}^N e^{-\phi^A_j}\abs{V_j}}^a,
 \end{equation*}
 with $(\hat{Z}_B)^b$ and $(\hat{Z}_C)^c$ defined analogously. $\hat{\pi}(a,b,c)$ denotes the equilibrium probability of having $(a,b,c)$ particles of species \textrm{A}, \textrm{B}, and \textrm{C} respectively. Substituting~\eqref{eq:Fss} into~\eqref{eq:correctDBMultPartF}, we can solve for $\km_{i j k}$ in terms of $\kp_{i j k}$ (or vice-versa)
\begin{equation}
	\km_{i j k} = \frac{ab}{c+1}\frac{\hat{\pi}(a,b,c)}{\hat{\pi}(a-1,b-1,c+1)}\frac{\hat{Z}_C}{\hat{Z}_A\hat{Z}_B}\frac{\abs{V_{ij}}}{\abs{V_{k}}}e^{-\phi^A_{i}-\phi^B_{j}+\phi^C_{k}}\kp_{i j k}.
	\label{eq:MultPartDBkmkppartial}
\end{equation}
To simplify the coefficient in the preceding relationship, we sum \eqref{eq:DBMultPartF} over all spatial locations (i.e. over all $\bi^a$, $\bj^b$, and $\bk^c$'s) to obtain
\begin{equation}
	\begin{aligned}
		0 = &-ab\kp\hat{\pi}(a, b, c) + (a+1)(b+1)\kp\hat{\pi}(a+1, b+1, c-1)\\
		&-c\km\hat{\pi}(a, b, c) + (c+1)\km\hat{\pi}(a-1, b-1, c+1).
	\end{aligned}
	\label{eq:MultPartDBWM}
\end{equation}
Here
\begin{equation*}
	\kp = \frac{1}{\hat{Z}_A\hat{Z}_B}\brac{\sum_{i = 1}^N\sum_{j = 1}^N\kpij\abs{V_{ij}}e^{-\phi^A_i-\phi^B_j}},
	\label{eq:MultPartkp}
\end{equation*}
and
\begin{equation*}
	\km = \frac{1}{\hat{Z}_C}\brac{\sum_{k = 1}^N\kmk\abs{V_{k}}e^{-\phi^C_k}},
	\label{eq:MultPartkm}
\end{equation*}
correspond to the average association and dissociation rates with respect to the Gibbs-Boltzmann distribution for a minimal set of substrates of each reaction. We define the non-dimensional dissociation constant $\Kd$ as their ratio,
\begin{equation*}
	\Kd = \frac{\km}{\kp}.
\end{equation*}
As we expect the equilibrium probabilities, $\hat{\pi}$, to also satisfy detailed balance, we should have that
\begin{equation*}
	\kp ab \, \hat{\pi}(a,b,c) = \km (c+1) \hat{\pi}(a-1,b-1,c+1)
\end{equation*}
and hence
\begin{equation*}
	\Kd = \frac{ab\hat{\pi}(a,b,c)}{(c+1)\hat{\pi}(a-1,b-1,c+1)}.
\end{equation*}
This is consistent with our choice that $\Kd = \hat{\pi}_{\textrm{A}\textrm{B}} / \hat{\pi}_{\textrm{C}}$ in the two-particle case studied in the main text (i.e. $a = 1$, $b = 1$, and $c = 0$).
\eqref{eq:MultPartDBkmkppartial} then simplifies to
\begin{equation*}
	\km_{i j k} = \frac{\Kd\ZhatC}{\ZhatA\ZhatB}\frac{\abs{V_{ij}}}{\abs{V_{k}}}e^{-\phi^A_{i}-\phi^B_{j}+\phi^C_{k}}\kp_{i j k}.
\end{equation*}
Notice this is identical to~\eqref{eq:discDBkmfromkp}. We therefore conclude that whenever the two particle association rate and the one-particle dissociation rate satisfy the detailed balance relation \eqref{eq:discDBkmfromkp}, the corresponding multiparticle transition rates will satisfy the (equivalent) detailed balance relation~\eqref{eq:MultPartDBkmkppartial}. That is, choosing $\kp_{i j k}$ and $\km_{i j k}$ to ensure detailed balance in the two particle case is sufficient to ensure consistency of the reactive transition rates with detailed balance holding in the general multiparticle case.

\end{appendix}

\end{document}